\documentclass[EJP]{ejpecp} 

\usepackage{array}
\usepackage{graphics}

\AUTHORS{%
  Nadia~Belaribi\footnote{Universit\'e Paris 13 and ENSTA ParisTech, France.
    \EMAIL{belaribi@math.univ-paris13.fr}}
  \and 
  Francesco~Russo\footnote{ENSTA ParisTech,
UMA,  828, Boulevard des Mar\'echaux,
F-91120 Palaiseau, France. 
  \EMAIL{francesco.russo@ensta-paristech.fr}} }

\SHORTTITLE{Fokker-Planck and Fast Diffusion}

\TITLE{Uniqueness for Fokker-Planck equations with measurable coefficients and applications to the fast diffusion equation.} 

\KEYWORDS{Fokker-Planck; fast diffusion; probabilistic representation;
non-linear diffusion; stochastic particle algorithm} 

\AMSSUBJ{MSC 2010:
 60H30; 60G44; 60J60; 60H07; 35C99; 35K10; 35K55; 35K65; 65C05; 65C35 } 

\SUBMITTED{November 28th, 2011} 

\VOLUME{0}
\YEAR{2012}
\PAPERNUM{0}
\DOI{vVOL-PID}

\ABSTRACT{The object of this paper is the uniqueness
for a $d$-dimensional Fokker-Planck type equation with
 inhomogeneous (possibly degenerated)
measurable  not necessarily bounded
coefficients. We provide an application to the
 probabilistic representation of  the so-called Barenblatt's solution
of the  fast diffusion equation which is
the partial differential equation
 $\partial_t u = \partial^2_{xx} u^m$ with $m\in]0,1[$.
Together with the mentioned Fokker-Planck equation,
we make use of   small
 time  density  estimates uniformly  with respect to the initial condition.}

\newtheorem*{hypothesis*}{Hypothesis}
\newcommand{\R}{\mathbb R}

\begin{document}



\section{Introduction}

The present paper is divided into three parts.
\begin{description}
  \item i) A uniqueness result
on a Fokker-Planck type equation with measurable non-negative
(possibly degenerated) multidimensional
  unbounded coefficients.
  \item ii) An  application to the probabilistic representation
 of a fast diffusion equation.
  \item iii) Some small time density estimates uniformly with respect to the initial condition.
\end{description}

In the whole paper $T > 0$
will stand for a fixed final time. In a one  dimension space, the  Fokker-Planck equation is of the type
\begin{equation}\label{FokkerPlank}
    \left\{
 \begin{array}{ccl}
  \partial_tu(t,x)&=& \partial_{xx}^2(a(t,x)u(t,x))-\partial_x(b(t,x)u(t,x)),\ \ t\in]0,T], \ x\in\R, \\
        u(0,\cdot)&=&\mu(dx), \\
\end{array}
\right.
\end{equation}
where  $a,b:[0,T]\times\R\rightarrow\R$ are measurable locally bounded coefficients and $\mu$ is a finite real Borel measure.
The Fokker-Planck equation for measures is a widely
studied subject in the literature whether in finite or
infinite dimension. Recent work in the case of time-dependent coefficients with some minimal regularity
was done by \cite{stannat,  figalli, Zhang} in the case $d \geq 1$. In infinite dimension
some interesting work was produced by \cite{BogaDaprato}.

In this paper we concentrate on
the case of measurable (possibly) degenerate coefficients.
Our interest is devoted to the irregularity of the diffusion
coefficient, so
 we will set $b = 0$.
A first result in that direction was produced in
\cite{BRR1} where $a$ was  bounded, possibly degenerated,
and the difference of two solutions was supposed to be in
$L^2([\kappa,T]\times\mathbb{R})$, for every $\kappa>0$ (ASSUMPTION (A)). This result was applied to  study
the probabilistic representation of a porous media type
equation with irregular coefficients. We will later come back
 to this point.
We remark that it is not possible to obtain uniqueness
without ASSUMPTION (A). In particular  \cite[Remark 3.11]{BRR1}
provides two measure-valued solutions
when $a$ is time-homogeneous, continuous,  with
$\frac{1}{a}$ integrable in a neighborhood of zero.

One natural question is about what happens when $a$ is not bounded
and $x \in \R^d$.
A partial answer to this question is given in Theorem \ref{TheoBRREtendu} which is probably
the most important result of the paper;
 it is a generalization of   \cite[Theorem 3.8]{BRR1} where
 the inhomogeneous function
$a$ was   bounded.
Theorem \ref{TheoBRREtendu}
 handles the multidimensional case and it
allows $a$ to be unbounded.

An application of Theorem \ref{TheoBRREtendu}  concerns
 the  parabolic problem:
\begin{equation}\label{PME}
    \left\{
 \begin{array}{ccl}
  \partial_tu(t,x)&=& \partial_{xx}^2(u^m(t,x)),\ \ t\in]0,T],\ x \in \R,\\
        u(0,\cdot)&=&\delta_{\scriptscriptstyle{0}}, \\
\end{array}
\right.
\end{equation}
where $\delta_{\scriptscriptstyle{0}}$ is the Dirac measure
at zero and $u^m$ denotes $u |u|^{m-1}$. It is well known that, for $m>1$, there exists an exact solution to
\eqref{PME}, the so-called \emph{Barenblatt's density}, see \cite{Barenb}.
Its explicit formula is recalled for instance in  \cite[Chapter 4]{LivreVasqz} and more precisely in  \cite[Section 6.1]{BCR1}. Equation \eqref{PME} is the {\it classical} porous medium equation.

In this paper, we  focus on \eqref{PME} when
$m\in ]0,1[$: the \emph{fast diffusion equation}. In fact, an analogous
 Barenblatt type solution also exists in this case, see \cite[Chapter 4]{LivreVasqz} and references therein; it is given by
the  expression
\begin{equation}\label{FormExpDirac}
    \mathcal{U}(t,x)=t^{-\alpha}{\left(D+\tilde{k}|x|^2t^{-2\alpha}\right)}^{-\frac{1}{1-m}},
\end{equation}
where
\begin{equation}\label{Constants}
\alpha=\frac{\displaystyle{1}}{\displaystyle{m+1}},~\tilde{k}=\frac{\displaystyle{1-m}}{\displaystyle{2(m+1)m}},~
D=\left(\frac{\displaystyle{I}}{\displaystyle{\sqrt{\tilde{k}}}}\right)^{\frac{2(1-m)}{m+1}},~I=
\int_{-\frac{\pi}{2}}^{\frac{\pi}{2}}
\left[\cos(x)\right]^{\frac{2m}{1-m}}dx.
\end{equation}

Equation \eqref{PME} is a particular case of the so-called
generalized porous media type equation
\begin{equation}\label{EDP}
 \left\{
 \begin{array}{ccl}
  \partial_tu(t,x)&=& \partial_{xx}^2 \beta\left(u(t,x)\right),\ \
t\in]0,T],\\
        u(0,x)&=&u_0(dx), \ \  \ x \in \mathbb{R},\\
\end{array}
\right.
\end{equation}
where $\beta: \mathbb{R} \rightarrow \mathbb{R}$ is a monotone
non-decreasing function
such that $\beta(0)=0$ and $u_0$ is a finite measure.  When $\beta(u) = u^m, \  m \in ]0,1[$  and   $u_0=\delta_0$,   two difficulties arise: first,
 the coefficient $\beta$ is of singular type since
it is not locally Lipschitz, second,  the initial condition
is a measure. Another type of singular coefficient
is $\beta(u) = H(u-u_c) u$, where $H$ is a Heaviside function and $u_c>0$ is some critical value, see e.g. \cite{BRR2}. Problem  \eqref{PME}  with $m\in ]0,1[$ was
 studied by several  authors. For a bounded integrable
function as initial condition, the equation in \eqref{PME}  is well-stated
in the sense of distributions,
as a by product of the classical papers  \cite{BenCran_79,
BenCran_81} on \eqref{EDP} with general monotonous coefficient
$\beta$. When the  initial data is locally integrable,   existence was proved by \cite{HerPier}.  \cite{BrezFried} extended the validity of this result when $u_0$ is a finite Radon measure in a bounded
domain, \cite{Pierre_87} established existence when $u_0$ is a locally finite measure in
 the whole space.  The Barenblatt's solution is an
 {\it extended continuous solution} as defined
in \cite{ChassVaz_02, ChassVaz_06};  \cite[Theorem 5.2]{ChassVaz_06} showed uniqueness in that class. \cite[Theorem 3.6]{TEEMU.L}
   showed  existence in a bounded domain of solutions to the fast diffusion equation
 perturbed by a right-hand side source term, being a general finite and positive Borel measure.
 As far as we know,  there is no
uniqueness argument in the literature whenever the initial
condition is a  finite measure in the general sense of distributions. Among  recent contributions,    \cite{Tosc} investigated the
large time behavior of solutions to \eqref{PME}.

The present paper provides the probabilistic representation
of the (Barenblatt's)  solution of \eqref{PME}
and exploits this fact in order to approach it via a
Monte Carlo simulation with an
 $L^2$ error  around $10^{-3}$.
 We make use of  the probabilistic procedure
 developed in   \cite[Section 4]{BCR1}  and we compare it to
 the exact form of the solution  $\mathcal{U}$  of \eqref{PME} which
 is given by the explicit formulae \eqref{FormExpDirac}-\eqref{Constants}. The target of \cite{BCR1} was the case $\beta(u) = H(u-u_c)u$; in that paper
  those techniques were compared  with a deterministic numerical analysis recently developed in \cite{Cavalli&Co} which was very performing
in that target case. At this stage, the implementation
of the same deterministic method for the fast diffusion equation
does not give satisfying results; this constitutes
a further justification for the probabilistic representation.

We define
\begin{equation*}
    \Phi(u)=|u|^{\frac{m-1}{2}}, \ \ u\in \mathbb{R}, \ \ m\in]0,1[.
\end{equation*}
The probabilistic
representation of $\mathcal{U}$ consists in finding a suitable
stochastic process $Y$
such that the law of $Y_t$ has  $\mathcal{U}(t,\cdot)$
as density.
$Y$ will be a (weak) solution of the
non-linear SDE
\begin{equation}\label{NLSDE}
    \left\{
 \begin{array}{ccl}
   Y_t&=&\int\limits_{0}^t \sqrt{2}\Phi(\mathcal{U}(s,Y_s))dW_s,\\
       \mathcal{U}(t,\cdot)&=& \mbox{Law  density of}
\ Y_t,\ \ \forall~t\in]0,T],\\
 \end{array}
\right.
\end{equation}
where  $W$ is a Brownian motion on some suitable filtered
probability space $(\Omega, \mathcal{F},(\mathcal{F}_t)_{t\geq
0},{P})$.

 To the best of our knowledge, the first author who considered a
probabilistic representation of a solution of \eqref{EDP}
was H. P. Jr. McKean (\cite{McKean}), particularly in relation with the
so-called \emph{propagation of
chaos}. In his case  $\beta$ was smooth, but the equation
also included a first order coefficient.
From then on, literature steadily grew and nowadays there is a vast amount of
contributions to the subject,  especially when the non-linearity is in the
first order part, as e.g. in Burgers' equation. We refer the reader to the
excellent survey papers \cite{Sznit} and \cite{graham}. A probabilistic interpretation of \eqref{EDP} when
$\beta(u)=u.|u|^{m-1},~m>1$, was provided for
instance in \cite{BCRV}. Recent developments related  to chaos propagation when $\beta(u)= u^2$
 and $\beta(u) = u^m,
  m > 1$ were proposed in
 \cite{Philipow} and \cite{Philipow_Figal}.
 The probabilistic representation in the case of possibly
discontinuous $\beta$ was treated in \cite{BRR1} when
$\beta$ is non-degenerate and in \cite{BRR2} when
$\beta$ is degenerate; the latter case includes the
case $\beta(u) = H(u-u_c) u$.

As a preamble to the probabilistic representation we
make a simple, yet crucial  observation.
  Let $W$ be a standard Brownian motion.

\begin{proposition}\label{PropI2}

Let $\beta:\mathbb{R}\rightarrow \mathbb{R}$  such that $\beta(u)=\Phi^2(u).u$, $\Phi:\mathbb{R}\rightarrow \mathbb{R}_+$ and $u_0$ be a probability
real measure.

Let $Y$ be a solution to the problem
\begin{equation}\label{CorI2_NLSDE}
    \left\{
 \begin{array}{ccl}
   Y_t&=&Y_0+\int\limits_{0}^t \sqrt{2}\Phi({u}(s,Y_s))dW_s,\\
       {u}(t,\cdot)&=& \mbox{Law  density of}
\ Y_t,\ \ \forall~t\in]0,T],\\
u(0,\cdot)&=&u_0(dx).\\
 \end{array}
\right.
\end{equation}
Then $u:[0,T]\times \mathbb{R}\rightarrow \mathbb{R}$ is solution to \eqref{EDP}.
\end{proposition}
Proof of the above result  is based on the following lemma.

\begin{lemma}\label{LemmaI1}
Let $a:[0,T]\times \mathbb{R}\rightarrow \mathbb{R}_+$ be measurable. Let $(Y_t)$ be a process which solves the SDE
\begin{equation*}
    Y_t=Y_0+\int\limits_0^t\sqrt{2a(s,Y_s)} dW_s, \ \ \ t\in[0,T].
\end{equation*}
Consider the function $t\mapsto \rho(t,\cdot)$ from $[0,T]$ to the space of finite real measures $\mathcal{M}(\mathbb{R})$, defined as $\rho(t,\cdot)$ being the law of $Y_t$. Then $\rho$ is a solution, in the sense of distributions (see \eqref{PDE_Th3.8}), of
\begin{equation}\label{PDE_a}
    \left\{
 \begin{array}{ccl}
  \partial_tu&=& \partial_{xx}^2(a u),\ t\in]0,T], \\
        u(0,\cdot)&=& \mbox{Law of $Y_0$}.
\end{array}
\right.
\end{equation}
\end{lemma}
\begin{proof}[Proof of Lemma \ref{LemmaI1}]

This is a classical result,  see for instance \cite[Chapter 4]{StrVarBook}.
 The proof is based on an application of It\^o's formula to
 $\varphi(Y_t)$, $\varphi \in \mathcal{S}(\mathbb{R})$.
\end{proof}

\begin{proof}[Proof of Proposition \ref{PropI2}]
We set $a(s,y)=\Phi^2(u(s,y))$. We apply Lemma \ref{LemmaI1}
setting $\rho(t,dy)=u(t,y)dy$, $t\in]0,T]$,  and $\rho(0,\cdot)=u_0$.
\end{proof}

When $u_0$ is the Dirac measure at zero and $\beta(u) = u^m$, with $m \in ]\frac{3}{5}, 1[$, Theorem \ref{Theorem(B)}
states the converse of Proposition \ref{PropI2},
providing a process $Y$ which is the unique (weak) solution of
\eqref{NLSDE}. The first step  consists in reducing  the proof of that Theorem
to the proof of Proposition \ref{Proposition(A)}  where
the Dirac measure, as initial condition of \eqref{PME},
 is replaced by the function $ \mathcal{U}(\kappa, \cdot),\  0 < \kappa \leq T$.
This corresponds to the shifted
Barenblatt's solution along a time $\kappa$, which will be
denoted by $ \overline{\mathcal{U}}$.
Also, in this case  Proposition \ref{Proposition(A)} provides
  an unique strong solution of the corresponding non-linear
SDE. That reduction is possible through a weak convergence argument
of the solutions given by Proposition \ref{Proposition(A)}
when $\kappa \rightarrow 0$. The idea of the proof of Proposition \ref{Proposition(A)}
is the following.
Let $W$ be a standard Brownian motion and
$\overline{Y_0}$  be  a r.v. distributed as $ \mathcal{U}(\kappa, \cdot)$; since
 $\Phi(\overline{\mathcal{U}})$ is Lipschitz, the
SDE  \[ \overline{Y_t} =  \overline{Y_0} + \int_0^t
\Phi(\overline{\mathcal{U}}(s, \overline{Y_s})) dW_s,\ t\in ]0,T], \]
admits a unique strong solution. The marginal laws of $(\overline{Y_t})$
and $ \overline{\mathcal{U}}$ can be shown to be both solutions
to \eqref{PDE_a} for $a(s,y) = (\overline{\mathcal{U}}(s,y))^{m-1}$;
that $a$ will be denoted in the sequel   by $\bar{a}$.
The leading argument of the proof is carried by
 Theorem \ref{TheoBRREtendu} which states uniqueness for
measure valued solutions of the  Fokker-Planck type  PDE \eqref{PDE_a}
 under some
 {\bf Hypothesis(B)}.
More precisely, to conclude that the marginal laws of  $ (\overline{Y_t})$
and  $\overline{\mathcal{U}}$ coincide via Theorem  \ref{TheoBRREtendu},
we show that they  both
verify the  so-called {\bf Hypothesis(B2)}.
In order to prove that for $\overline{\mathcal{U}}$,   we will make use
of Lemma \ref{Lemma(A1)}.
The verification of {\bf Hypothesis(B2)}  for the  marginal laws
of $\overline{Y}$ is more involved. It makes use of a small time
(uniformly with respect to the initial condition)
  upper bound for the density of an inhomogeneous diffusion flow
with linear growth  (unbounded) smooth coefficients,
even though the diffusion term is non-degenerate and all the derivatives are bounded.
This is the object of Proposition \ref{Prop_UpperBound_Dens},
the proof of which is based on an application of Malliavin calculus.
 In our opinion  this result alone is  of interest as
we were not able to find it in the literature.\\
When the paper was practically finished we  discovered
an interesting recent result of M. Pierre, presented in
\cite[Chapter 6]{profeta},  obtained independently. This result holds in dimension $1$
when the   coefficients  are locally bounded, non-degenerate and the initial condition has a first moment.
In this  case, the hypothesis of type (B) is not needed. In particular it allows  one to establish Proposition \ref{Proposition(A)},  but not Theorem \ref{Theorem(B)}  where the coefficients are not locally bounded on $[0,T]\times \R$.

The paper is organized as follows. Section \ref{Prelim}  is devoted to
basic notations.
Section \ref{LinPDE} is concentrated on Theorem \ref{TheoBRREtendu}
which concerns uniqueness for the deterministic, time inhomogeneous,
Fokker-Planck type equation.
Section \ref{Basic_U}  presents some properties of the Barenblatt's solution
 $\mathcal{U}$ to \eqref{PME}.
The probabilistic representation of  $\mathcal{U}$
 is treated in
Section \ref{ProbRep}.
Proposition \ref{Prop_UpperBound_Dens}
performs  small time density estimates
for time-inhomogeneous
diffusions, the  proof of which is located in the Appendix.
Finally, Section \ref{NumExp} is devoted to numerical experiments.



\section{Preliminaries}\label{Prelim}
We start with some basic analytical framework. In the whole paper $d$ will be a strictly positive integer.
If $f:\mathbb{R}^d\rightarrow \mathbb{R}$ is a bounded function we
will denote  $\|f\|_{\infty}=\sup\limits_{x\in \mathbb{R}^d}|f(x)|$.
By $S(\mathbb{R}^d) $ we denote the space of rapidly decreasing
infinitely differentiable functions $\varphi: \mathbb{R}^d\rightarrow
\mathbb{R}$,  by $S^{\prime}(\mathbb{R}^d)  $ its dual (the
space of tempered distributions). We denote by $\mathcal{M}(\mathbb{R}^d)$  the set of finite Borel  measures on $\R^d$. If $x\in \R^d$, $|x|$ will denote the usual Euclidean norm.

For  $\varepsilon > 0$,
let $K_\varepsilon$ be the Green's function of $\varepsilon - \Delta$,
that is the kernel of the operator  $(\varepsilon - \Delta)^{-1}:
L^2(\mathbb{R}^d) \rightarrow H^2(\mathbb{R}^d)\subset L^2(\mathbb{R}^d)$.
 In particular, for all $\varphi \in L^2(\mathbb{R}^d)$, we have
\begin{equation}
B_\varepsilon\varphi :=
   (\varepsilon  -\Delta)^{-1}\varphi  (x)=\int_{\mathbb{R}}%
K_\varepsilon \left(  x-y\right)  \varphi(y)dy.\label{eq kernelbis}%
\end{equation}
For more information about the corresponding analysis,
 the reader can consult \cite{Stein}.
If $\varphi \in C^2 (\mathbb{R}^d) \bigcap S^{\prime}(\mathbb{R}^d)$,
then $(\varepsilon  - \Delta) \varphi$
 coincides with the classical associated PDE operator evaluated at $\varphi$.

\begin{definition}\label{DefNonDeg}

We will say that a function $\psi:[0,T]\times \mathbb{R}\rightarrow\mathbb{R}$ is {non-degenerate} if there is
a constant $c_{\scriptscriptstyle{0}}>0$  such that $\psi\geq
c_{\scriptscriptstyle{0}}$.
\end{definition}

\begin{definition}\label{LinGrowth}
We will say that a function $\psi:[0,T]\times \mathbb{R}\rightarrow\mathbb{R}$ has linear growth (with respect to the second variable)  if there is
a constant $C$  such that $|\psi(\cdot,x)|\leq C(1+|x|)$, $x\in \R$.
\end{definition}

\begin{definition}\label{def: sol_distrib}
Let $a:[0,T]\times \R^d\rightarrow \R_+$ be a Borel function,  $z^0\in \mathcal{M}(\R^d)$. A (weakly measurable) function $z:[0,T]\rightarrow \mathcal{M}(\R^d)$ is said to be a solution in the sense of distributions of \[\partial_t
z=\Delta(az)\]
with initial condition $z(0,\cdot)=z^0$  if,  for every $t\in[0,T]$ and $\phi\in \mathcal{S}(\R)$,  we have
\begin{equation}\label{PDE_Th3.8}
\int_{\mathbb{R}^d}\phi(x)z(t,dx)=\int_{\mathbb{R}^d}\phi(x)z^0(dx)+\int_0^t
ds \int_{\mathbb{R}^d}\Delta\phi(x)a(s,x)z(s,dx).
\end{equation}
\end{definition}

\section{Uniqueness for the Fokker-Planck equation}\label{LinPDE}

We now state  the main result of the paper
which concerns uniqueness for the Fokker-Planck type
equation with measurable, time-dependent,
(possibly degenerated and unbounded) coefficients. It generalizes  \cite[Theorem 3.8]{BRR1}
where the coefficients were bounded and one-dimensional.

The theorem below holds with two classes of hypotheses: \textbf{(B1)}, operating in the multidimensional case,  and \textbf{(B2)},
more specifically in the one-dimensional case.

\begin{theorem}\label{TheoBRREtendu}

Let $a$ be a Borel nonnegative  function on $[0,T]\times
\mathbb{R}^d$. Let $z_i:[0,T]\rightarrow \mathcal{M}(\mathbb{R}^d)$, $i=1,2$, be
continuous with respect to the weak topology on finite measures on
$\mathcal{M}(\mathbb{R}^d)$. Let $z^0$ be an element of $\mathcal{M}(\mathbb{R}^d)$. Suppose that
both $z_1$ and $z_2$ solve the problem $\partial_t
z=\Delta(az)$ in the sense of distributions with initial
condition $z(0,\cdot)=z^0$.

Then $z:=(z_1-z_2)(t,\cdot)$ is identically zero for every $t$
under the following requirement.
\begin{hypothesis*}[\textbf{B}]
There is $\tilde{z}\in L_{{loc}}^1([0,T]\times \mathbb{R}^d)$  such that $z(t,\cdot)$ admits $\tilde{z}(t,\cdot)$ as density for almost all $t\in[0,T]$;   $\tilde{z}$ will still be denoted  by $z$.
Moreover, either  \textbf{(B1)}
or \textbf{(B2)} below is fulfilled.\\
\[\text{\textbf{(B1)}}\quad \text{(i)}\int\limits_{[0,T]\times\mathbb{R}^d} \!|z(t,x)|^2 \,dt \,dx<+\infty, \quad \text{(ii)} \int\limits_{[0,T]\times\mathbb{R}^d}|az|^2(t,x)dt dx<+\infty.\]
\textbf{(B2)} We suppose $d=1$. For every $t_{\scriptscriptstyle{0}}>0$, we have
\[\text{(i)}\int\limits_{[t_{\scriptscriptstyle{0}},T]\times\mathbb{R}} \!|z(t,x)|^2
\,dt \,dx<+\infty, \quad \text{(ii)}\int\limits_{[0,T]\times\mathbb{R}} \!|az|(t,x) \,dt \,dx<+\infty, \quad \text{(iii)}\int\limits_{[t_{\scriptscriptstyle{0}},T]\times\mathbb{R}}\!|az|^2(t,x)\,dt \,dx<+\infty.\]
\end{hypothesis*}
\end{theorem}

\begin{remark}\label{R_Weak}
\({}\)
The weak continuity of $z(t,\cdot)$ and  \cite[Remark 3.10]{BRR1} imply that
$\sup\limits_{t\in[0,T]}\|z(t,\cdot)\|_{{var}}<+\infty$,  where $\|\cdot\|_{{var}}$ denotes the total variation. In particular  $\sup\limits_{0<t\leq T} \int_{\mathbb{R}^d}|z(t,x)|dx <+\infty$.
\end{remark}
\begin{remark}\label{RemarkConsqMultiDim}
\begin{enumerate}
\item If $a$ is bounded  then the first item of Hypothesis(B1) implies the second one.
\item If $a$ is non-degenerated, assumption (ii) of Hypothesis(B1) implies assumption (i).
\end{enumerate}
\end{remark}
\begin{remark}\label{RemarkConsq}

Let $d=1$.

\begin{enumerate}
    \item  If $a$ is non-degenerate, the third assumption of Hypothesis(B2) implies the first one.

    \item   If  $ z(t,x)\in L^{\infty}([t_{\scriptscriptstyle{0}},T]\times
    \mathbb{R})$  then the first item of Hypothesis(B2) is always verified.

    \item  If $a$ is bounded  then  assumption (ii) of Hypothesis(B2) is always verified by Remark \ref{R_Weak}; the first item of Hypothesis(B2) implies the third one. So  Theorem \ref{TheoBRREtendu}  is a
    strict generalization of  \cite[Theorem 3.8]{BRR1}.

    \item   Let  $(z(t,\cdot), t\in[0,T])$ be the marginal law densities  of a
    stochastic process $Y$ solving
    \begin{equation*}
        Y_t=Y_{\scriptscriptstyle{0}}+\int_0^t \sqrt{2a(s,Y_s)}dW_s,
    \end{equation*}
    with  $Y_{\scriptscriptstyle{0}}$ distributed as $z^0$  such that
    $\int\limits_{\mathbb{R}}|x|^2z^0(dx)<+\infty$.

If $\sqrt{a}$ has linear growth, it is well known that $\sup\limits_{t\leq T}\mathds{E}(|Y_t|^2)<+ \infty$; so
\begin{equation*}
\int\limits_{[0,T]\times\mathbb{R}} \!|a(s,x)z(s,x)| \,ds \,dx=\mathds{E}\left[\int\limits_0^Ta(s,Y_s)\,ds\right]<+\infty.
\end{equation*}
Therefore  assumption (ii) in Hypothesis(B2) is always fulfilled.
\end{enumerate}

\end{remark}

\begin{proof}[Proof of Theorem \ref{TheoBRREtendu}]

 Let $z_1$, $z_2$ be two solutions of \eqref{PDE_Th3.8}; we set
$z:=z_1-z_2$. We evaluate,  for every $t\in[0,T]$, the quantity
\begin{equation*}
    g_{\varepsilon}(t)=\|z(t,\cdot)\|^2_{-1,\varepsilon},
\end{equation*}
where
$\|f\|_{-1,\varepsilon}=\|(\varepsilon-\Delta)^{-\frac{1}{2}}f\|_{L^2}$.

Similarly to the first part of the proof of  \cite[Theorem 3.8]{BRR1}, assuming we can show that
\begin{equation}\label{Limite_t_zero}
    \lim_{\varepsilon\to 0}g_{\varepsilon}(t)=0,\ \ \ \forall t\in[0,T],
\end{equation}
we are able to prove that $z(t)\equiv 0$ for all $t\in]0,T]$.
We explain this fact.

Let $t\in]0,T]$. We recall the notation
$B_{\varepsilon}f=(\varepsilon-\Delta)^{-1}f$, if $f\in L^2(\mathbb{R}^d)$.
 Since $z(t,\cdot) \in L^2(\R^d)$  then $B_{\varepsilon}z(t,\cdot)\in H^2(\R^d)$
 and so $\nabla B_{\varepsilon}z(t,\cdot) \in H^1(\mathbb{R}^d)^d \subset L^2(\mathbb{R}^d)^d$. This gives
\begin{eqnarray*}
  g_{\varepsilon}(t) &=& \int_{\mathbb{R}^d}B_{\varepsilon}z(t,x) z(t,x)dx =\varepsilon \int_{\mathbb{R}^d}(B_{\varepsilon}z(t,x))^2dx-\int_{\mathbb{R}^d}B_{\varepsilon}z(t,x)\Delta B_{\varepsilon}z(t,x)dx\\
  &=& \varepsilon \int_{\mathbb{R}^d}(B_{\varepsilon}z(t,x))^2dx+\int_{\mathbb{R}^d}|\nabla B_{\varepsilon}z(t,x)|^2dx.
\end{eqnarray*}
Since the two terms of the above  sum are non-negative, if \eqref{Limite_t_zero} holds,  then
$\sqrt{\varepsilon}B_{\varepsilon}z(t,\cdot) \rightarrow 0$  (resp. $|\nabla B_{\varepsilon}z(t,\cdot)| \rightarrow 0$) in $L^2(\mathbb{R}^d)$ (resp. in $L^2(\mathbb{R}^d)^d$). So, for all $t\in]0,T]$,  $z(t,\cdot)=\varepsilon B_{\varepsilon}z(t,\cdot)-\Delta B_{\varepsilon}z(t,\cdot) \rightarrow 0$,
in the sense of distributions, as $\varepsilon$ goes to zero. Therefore  $z\equiv 0$.

We proceed now with the proof of \eqref{Limite_t_zero}. We have the following identities  in the sense of
distributions:
\begin{eqnarray}
  z(t,\cdot) = \int_0^t\Delta(az)(s,\cdot)ds =\int_0^t(\Delta-\varepsilon)(az)(s,\cdot)ds+\varepsilon\int_0^t(az)(s,\cdot)ds,\label{Eq1}
\end{eqnarray}
which implies
\begin{eqnarray}
  B_{\varepsilon}z(t,\cdot) &=& -\int_0^t(az)(s,\cdot)ds+\varepsilon\int_0^tB_{\varepsilon}(az)(s,\cdot)ds.\label{Eq2}
\end{eqnarray}
Let $\delta>0$ and $(\phi_{\delta})$ a sequence of mollifiers
converging to the Dirac delta function at zero. We set $z_{\delta}(t,x)=\int_{\mathbb{R}^d}z(t,y)\phi_{\delta}(x-y)dy$, observing that $z_{\delta}\in (L^1\bigcap L^{\infty})([0,T]\times\mathbb{R}^d)$. Moreover, \eqref{Eq1}  gives
\begin{equation*}
    z_{\delta}(t,\cdot)=\int_0^t\Delta(az)_{\delta}(s,\cdot)ds.
\end{equation*}
 We suppose now Hypothesis(B1) (resp. (B2)). Let $t_{\scriptscriptstyle{0}}=0$ (resp. $t_{\scriptscriptstyle{0}}>0$). By assumption (B1)(ii) (resp. (B2)(iii)), we have  $\Delta(az)_{\delta}\in L^2([t_{\scriptscriptstyle{0}},T]\times\mathbb{R}^d)$. Thus, $z_{\delta}$ can be seen as a
function belonging to $C([t_{\scriptscriptstyle{0}},T];L^2(\mathbb{R}^d))$. Besides, identities \eqref{Eq1} and \eqref{Eq2} lead to
\begin{eqnarray}
  z_{\delta}(t,\cdot)
  &=&z_{\delta}(t_{\scriptscriptstyle{0}},\cdot)+\int_{t_{\scriptscriptstyle{0}}}^t(\Delta-\varepsilon)(az)_{\delta}(s,\cdot)ds+\varepsilon\int_{t_{\scriptscriptstyle{0}}}^t(az)_{\delta}(s,\cdot)ds,\label{Eq3}\\
  B_{\varepsilon}z_{\delta}(t,\cdot) &=&B_{\varepsilon}z_{\delta}(t_{\scriptscriptstyle{0}},\cdot) -\int_{t_{\scriptscriptstyle{0}}}^t(az)_{\delta}(s,\cdot)ds+\varepsilon\int_{t_{\scriptscriptstyle{0}}}^tB_{\varepsilon}(az)_{\delta}(s,\cdot)ds.\ \ \label{Eq4}
\end{eqnarray}
Proceeding through integration by parts with values in $L^2(\mathbb{R}^d)$, we get
\begin{align}
  \|z_{\delta}(t,\cdot)\|^2_{-1,\varepsilon}&-\|z_{\delta}(t_{\scriptscriptstyle{0}},\cdot)\|^2_{-1,\varepsilon}
   = -2\int_{t_{\scriptscriptstyle{0}}}^t ds <z_{\delta}(s,\cdot),(az)_{\delta}(s,\cdot)>_{L^2}\nonumber\\
   &\label{Eq5}\\
   &+ 2 \varepsilon\int_{t_{\scriptscriptstyle{0}}}^t ds <(az)_{\delta}(s,\cdot),B_{\varepsilon}z_{\delta}(s,\cdot)>_{L^2}.\nonumber
\end{align}
Then, letting $\delta$ go to zero, using assumptions (B1)(i)-(ii) (resp.  (B2)(i) and (B2)(iii)) and Cauchy-Schwarz inequality, we obtain
\begin{align}
  \|z(t,\cdot)\|^2_{-1,\varepsilon} &-\|z(t_{\scriptscriptstyle{0}},\cdot)\|^2_{-1,\varepsilon}=- 2\int_{t_{\scriptscriptstyle{0}}}^t ds \int_{\mathbb{R}^d}a(s,x)|z|^2(s,x)dx \nonumber\\
  &\label{Eq6}\\
   &+ 2 \varepsilon\int_{t_{\scriptscriptstyle{0}}}^t ds
   <(az)(s,\cdot),B_{\varepsilon}z(s,\cdot)>_{L^2}.\nonumber
\end{align}

At this stage of the proof, we assume  that \textbf{Hypothesis(B1)} is satisfied. Since $t_{\scriptscriptstyle{0}}=0$, we have $z(t_{\scriptscriptstyle{0}},\cdot)=0$. Using the inequality $c_{\scriptscriptstyle{1}}c_{\scriptscriptstyle{2}}\leq \frac{c_{\scriptscriptstyle{1}}^2+c_{\scriptscriptstyle{2}}^2}{2}$, $c_{\scriptscriptstyle{1}}$, $c_{\scriptscriptstyle{2}} \in \R$  and Cauchy-Schwarz, \eqref{Eq6} implies
\begin{eqnarray}
  \|z(t,\cdot)\|^2_{-1,\varepsilon} &\leq& -2\int_0^tds\int_{\mathbb{R}^d}(a|z|^{\scriptscriptstyle{2}})(s,x) dx + \varepsilon \int_0^tds \|az(s,\cdot)\|_{L^2}^2 + \varepsilon\int_0^tds \|B_{\varepsilon}z(s,\cdot)\|_{L^2}^2\nonumber\\
  &\leq & \varepsilon \int_0^tds \|az(s,\cdot)\|_{L^2}^2+ \int_0^tds \|z(s,\cdot)\|^2_{-1,\varepsilon},\label{Eq6Bis}
\end{eqnarray}
because for $f=z(s,\cdot)$, we have
\begin{equation*}
\varepsilon \|B_{\varepsilon}f\|_{L^2}^2=\varepsilon\int_{\mathbb{R}^d}\frac{\displaystyle{(\mathcal{F}(f))^2(\xi)}}{\displaystyle{(\varepsilon+|\xi|^2)^2}}d\xi
\leq  \int_{\mathbb{R}^d}\frac{\displaystyle{(\mathcal{F}(f))^2(\xi)}}{\displaystyle{\varepsilon+|\xi|^2}}d\xi=\|f\|^2_{-1,\varepsilon}.
\end{equation*}
We observe that the first integral of the right-hand side of \eqref{Eq6Bis} is finite by assumption (B1)(ii). Gronwall's lemma,  applied to \eqref{Eq6Bis}, gives
\begin{equation*}
    \|z(t,\cdot)\|^2_{-1,\varepsilon}\leq \varepsilon e^T \int_0^T ds \|az(s,\cdot)\|_{L^2}^2.
\end{equation*}
Letting $\varepsilon\rightarrow 0$, it follows that $\|z(t,\cdot)\|^2_{-1,\varepsilon}=0,\ \forall t\in[0,T]$. This concludes the first part of the proof.

We now suppose   that \textbf{Hypothesis(B2)} is satisfied, in particular $d=1$. By  \cite[Lemma 2.2]{BRR1}  we have
\begin{equation*}
    \sup_{x} 2 \varepsilon|B_{\varepsilon}z(s,x)|\leq
    \sqrt{\varepsilon}\|z(s,\cdot)\|_{var}.
\end{equation*}
Consequently  \eqref{Eq6} gives
\begin{equation}\label{Eq7}
  \|z(t,\cdot)\|^2_{-1,\varepsilon} -\|z(t_{\scriptscriptstyle{0}},\cdot)\|^2_{-1,\varepsilon}
   \leq  \sqrt{\varepsilon}\sup_{t\leq T}\|z(t,\cdot)\|_{var}\int\limits_{[{t_{\scriptscriptstyle{0}}},T]\times \mathbb{R}}|az|(s,x) ds
   dx.
\end{equation}

Besides, arguing like in the proof of \cite[Theorem 3.8]{BRR1}, we
obtain that
\begin{equation*}
\lim\limits_{t_{\scriptscriptstyle{0}}\to
0}\|z(t_{\scriptscriptstyle{0}},\cdot)\|^2_{-1,\varepsilon}=0.
\end{equation*}
We first let $t_{\scriptscriptstyle{0}}\to 0$ in \eqref{Eq7}, which implies
\begin{equation}\label{Eq7Bis}
   \|z(t,\cdot)\|^2_{-1,\varepsilon}
   \leq  \sqrt{\varepsilon}\sup_{t\leq T}\|z(t,\cdot)\|_{{var}}\int\limits_{[0,T]\times \mathbb{R}}|az|(s,x) ds;
\end{equation}
we remark that the right-hand side of \eqref{Eq7Bis} is finite by assumption (B2)(ii). Letting $\varepsilon$ go to zero,
the proof of \eqref{Limite_t_zero} is finally established.

\end{proof}

\section{Basic facts on the fast diffusion equation}\label{Basic_U}
We go on providing some properties of the Barenblatt's solution $\mathcal{U}$ to \eqref{PME} when $m\in]0,1[$ and given by \eqref{FormExpDirac}-\eqref{Constants}.
\begin{proposition}\label{Prop_U_pp}

\({}\)

\begin{description}
\item (i) $\mathcal{U}$ is a solution in the sense of distributions to \eqref{PME}. In particular, for every $\varphi \in C_0^{\infty}(\mathbb{R})$,   we have
\begin{equation}\label{DistribEq}
    \int\limits_{\mathbb{R}}\varphi(x)\mathcal{U}(t,x)dx=\varphi(0)+\int\limits_0^tds\int\limits_{\mathbb{R}}\mathcal{U}^m(s,x)\varphi''(x)dx.
\end{equation}
\item  (ii) $\int\limits_{\mathbb{R}}\mathcal{U}(t,x)dx=1,  \ \ \ \forall t>0$. 
In particular, for any $t > 0$, $\mathcal{U}(t,\cdot)$
is a probability density.
\item (iii) The Dirac measure $\delta_{\scriptscriptstyle{0}}$  is the initial trace of $\mathcal{U}$, in the sense that
\begin{equation}\label{M2}
    \int\limits_{\mathbb{R}}\gamma(x)\mathcal{U}(t,x)dx\to \gamma(0), \mbox{ as }\  t\to 0,
\end{equation}
for every $\gamma:\mathbb{R}\rightarrow \mathbb{R}$, continuous and bounded.
\end{description}
\end{proposition}
\begin{proof}[Proof of Proposition \ref{Prop_U_pp}]
(i)  This is a well known fact which can be established by inspection.

(ii)  For $M\geq 1$, we consider a sequence of smooth functions $(\varphi^M)$, such that
 \[\varphi^M(x)\left\{
                      \begin{array}{ll}
                        =0, & \hbox{if $|x|\geq M+1$;} \\
                        \leq 1, & \hbox{if $|x|\in[M,M+1]$;} \\
                        =1, & \hbox{if $|x|\leq M$.}
                      \end{array}
                    \right.\]
 By \eqref{DistribEq}  we have
 \begin{equation}\label{DistribEq1}
    \int\limits_{\mathbb{R}}\varphi^M(x)\mathcal{U}(t,x)dx=1+\int\limits_0^tds\int\limits_{\mathbb{R}}\mathcal{U}^m(s,x)(\varphi^M)''(x)dx.
\end{equation}
Letting $M\to +\infty$, by Lebesgue's dominated convergence theorem, the left-hand side of \eqref{DistribEq1} converges to $\int_{\mathbb{R}}\mathcal{U}(t,x)dx$. The integral on the right-hand side of \eqref{DistribEq1} is bounded by
\begin{equation*}
  C \int\limits_0^tds\int\limits_{M}^{M+1}\mathcal{U}^m(s,x)dx \leq  C \int\limits_0^t  s^{-\alpha m}\left(D+\tilde{k}M^2s^{-2\alpha}\right)^{\frac{-m}{1-m}}ds\leq\frac{\displaystyle{C}}{\displaystyle{(\tilde{k}M^2)^{\frac{m}{1-m}}}}\int\limits_0^T  s^{{\frac{m}{1-m}}}ds.
\end{equation*}
The last  integral  on the right  is finite as $\frac{\displaystyle{m}}{\displaystyle{1-m}}>0$, for every $m\in]0,1[$. Therefore  the integral in the right-hand side of \eqref{DistribEq1} goes to zero as $M\to +\infty$. This concludes the proof of the second item of Proposition \ref{Prop_U_pp}.

(iii) \eqref{M2} follows by elementary changes of variables.
\end{proof}

Note that  the second item  of Proposition \ref{Prop_U_pp} determines the explicit expression of
the constant $D$.\\

\begin{lemma}\label{Lemma(A1)}

\({}\)

\begin{description}
\item (i)  Suppose that  $\ \frac{\displaystyle{1}}{\displaystyle{3}}<m<1\ $. Then there is $p\geq 2$  and a constant $\mathcal{C}_p$ (depending on $T$) such that for $0\leq s<\ell\leq T$
\begin{equation}\label{L(A1).(i).1}
    \int\limits_{]s,\ell]\times \mathbb{R}}dtdx\left(\mathcal{U}(t,x)\right)^{\scriptscriptstyle{\frac{p(m-1)}{2}+1}}\leq \mathcal{C}_p(\ell-s).
\end{equation}
\item (ii) In particular, taking  $p=2$ in \eqref{L(A1).(i).1},  we get
\begin{equation}\label{L(A1).(i).2}
    \int\limits_{]0,T]\times \mathbb{R}}dtdx\left(\mathcal{U}(t,x)\right)^{m}<+\infty,
\end{equation}
again when $m$ belongs to $]\frac{1}{3},1[$.
\item (iii) If $\ \frac{\displaystyle{1}}{\displaystyle{5}}<m<1\ $,
\begin{equation}\label{L(A1).(i).3}
    \int\limits_{]0,T]\times \mathbb{R}}dtdx\left(\mathcal{U}(t,x)\right)^{2m}<+\infty .
\end{equation}
\item (iv) If $m$ belongs to
    $]\frac{\displaystyle{3}}{\displaystyle{5}},1[\ $, then
    \begin{equation}\label{L(A1).(ii).1}
        \forall \kappa>0, \ \ \int\limits_{\mathbb{R}}|x|^{\scriptscriptstyle{4}}\mathcal{U}(\kappa,x)dx<+\infty.
    \end{equation}
\end{description}

\end{lemma}

\begin{proof}[Proof of Lemma \ref{Lemma(A1)}]

\({}\)

(i) Using \eqref{FormExpDirac}, we have
\begin{equation*}
 \ \     \int\limits_{]s,\ell]\times
\mathbb{R}}{\left(\mathcal{U}(t,x)\right)}^{\scriptscriptstyle{\frac{p(m-1)}{2}+1}}
dtdx=\int\limits_{]s,\ell]\times
    \mathbb{R}}t^{\scriptscriptstyle{\frac{-\alpha p(m-1)}{2}-\alpha}}{\left(D+\tilde{k}|x|^2t^{-2\alpha}\right)}^{\scriptscriptstyle{\frac{p}{2}-\frac{1}{1-m}}}dtdx.
\end{equation*}
Then, setting
$y=t^{-\alpha}x\sqrt{\frac{\displaystyle{\tilde{k}}}{\displaystyle{D}}}$,
we get
\begin{align*}
     \int\limits_{]s,\ell]\times
\mathbb{R}}{\left(\mathcal{U}(t,x)\right)}^{\scriptscriptstyle{\frac{p(m-1)}{2}+1}}
dtdx&=\frac{\displaystyle{D^{\scriptscriptstyle{\frac{p+1}{2}-\frac{1}{1-m}}}}}{\displaystyle{\sqrt{\tilde{k}}}}
    \int\limits_s^l t^{\frac{p}{2}\alpha(1-m)}dt\int_{\mathbb{R}}(1+y^2)^{\scriptscriptstyle{\frac{p}{2}-\frac{1}{1-m}}}dy\\
    &\leq \frac{\displaystyle{D^{\scriptscriptstyle{\frac{p+1}{2}-\frac{1}{1-m}}}}}{\displaystyle{\sqrt{\tilde{k}}}}T^{\frac{p}{2}\alpha(1-m)}(\ell-s)\int_{\mathbb{R}}(1+y^2)^{\scriptscriptstyle{\frac{p}{2}-\frac{1}{1-m}}}dy.
\end{align*}
The last integral is finite if $(p+1)(1-m)<2$. This implies \eqref{L(A1).(i).1}.\\

(ii) is a particular case of (i)  and (iii) follows by similar arguments as for the proof of (i).

(iv) Now  we  assume that $m\in]\frac{\displaystyle{3}}{\displaystyle{5}},1[$. For $\kappa>0$  we have
\begin{equation}\label{etoile(2)}
\int\limits_{\mathbb{R}}|x|^{\scriptscriptstyle{4}}\mathcal{U}(\kappa,x)dx=
\frac{\displaystyle{D^{\scriptscriptstyle{\frac{3-5m}{2(1-m)}}}}}{\displaystyle{{\tilde{k}}^{5/2}}}
\kappa^{\scriptscriptstyle{4\alpha}}\int\limits_{\mathbb{R}}|y|^{\scriptscriptstyle{4}}(1+y^2)^{-\frac{{1}}{{1-m}}}dy,
\end{equation}
where  this last equality was obtained setting $y=\kappa^{-\alpha}x\sqrt{\frac{\displaystyle{\tilde{k}}}{\displaystyle{D}}}$. Clearly, since $m\in]\frac{\displaystyle{3}}{\displaystyle{5}},1[$, the integral in  the right-hand side of \eqref{etoile(2)} is finite. Therefore \eqref{L(A1).(ii).1} is fulfilled.
\end{proof}

Let $\kappa \in ]0,T]$. Given $u:[0,T]\times \mathbb{R}\rightarrow \mathbb{R}$  we associate
\begin{equation}\label{Trans_u}
    \overline{u}(t,x)=u(t+\kappa,x), \ \ \ (t,x)\in[0,T-\kappa]\times\mathbb{R}.
\end{equation}
In particular  we have
\begin{equation}\label{Trans_U}
    \overline{\mathcal{U}}(t,x)=\mathcal{U}(t+\kappa,x).
\end{equation}
Moreover, for every $x\in \mathbb{R}$, we denote
\begin{equation}\label{U_kappa}
    u_{\scriptscriptstyle{0,\kappa}}(x)=\mathcal{U}(\kappa,x).
\end{equation}
\begin{remark}
Function $\overline{\mathcal{U}}$ solves the problem
\begin{equation}\label{Ubar_PDE}
    \left\{
 \begin{array}{ccl}
  \partial_tu&=& \partial_{xx}^2(u^m),\\
  u(0,\cdot)&=&u_{0,\kappa}\ . \\
\end{array}
\right.
\end{equation}
\end{remark}

\section{The probabilistic representation of the fast diffusion equation}\label{ProbRep}
We are now interested in a non-linear stochastic differential
equation rendering the probabilistic representation related to
\eqref{PME} and given by \eqref{NLSDE}. Suppose for a moment that $Y_{\scriptscriptstyle{0}}$ is a random
variable distributed according to $\delta_{\scriptscriptstyle{0}}$, so  $Y_{\scriptscriptstyle{0}}=0$ a.s. We recall that,  if there exists a process $Y$ being a solution in law of
\eqref{NLSDE}, then Proposition \ref{PropI2} implies that $u$
solves \eqref{PME} in the sense of distributions.

In this subsection we shall prove existence and uniqueness of solutions in
law  for \eqref{NLSDE}. In this respect  we first state a tool, given by Proposition \ref{Prop_UpperBound_Dens} below,
concerning the existence of an upper bound for the marginal law densities of
the solution $Y$ of an inhomogeneous SDE with unbounded coefficients. This result has an independent interest.

\begin{proposition}\label{Prop_UpperBound_Dens}
Let $\sigma, \ b:[0,T]\times \mathbb{R}\rightarrow \mathbb{R}$ be
continuous (not necessarily bounded) functions  such that $\sigma(t,\cdot)$, $b(t,\cdot)$ are smooth with bounded
derivatives of orders greater or equal than one. $\sigma$ is supposed to be non-degenerate.

Let $x_{\scriptscriptstyle{0}} \in \mathbb{R}$ and
$Y_t=(Y_t^{x_{\scriptscriptstyle{0}}})_{t\in[0,T]}$  be  the
solution of
\begin{equation}\label{SDE_Upper_Bound}
    Y_t=x_{\scriptscriptstyle{0}}+\int_0^t \sigma(r,Y_r)dW_r+ \int_0^t b(r,Y_r)dr.
\end{equation}
Then, for every $s>0$,  the law of  $Y_s$  admits a density
denoted $ p_s(x_{\scriptscriptstyle{0}},\cdot)$.

Moreover, we have
\begin{equation}\label{DensUpperBound}
    p_s(x_{\scriptscriptstyle{0}},x) \leq \frac{\displaystyle{K}}{\displaystyle{\sqrt{s}}}\left(1+|x_{\scriptscriptstyle{0}}|^{4}\right), \ \forall (s,x)\in ]0,T]\times \mathbb{R},
\end{equation}
where  $K$ is a  constant which depends on $\|\sigma'\|_{\scriptscriptstyle{\infty}}$, $\|b'\|_{\scriptscriptstyle{\infty}}$ and $T$ but not on
$x_{\scriptscriptstyle{0}}$.
\end{proposition}

\begin{remark}
\begin{enumerate}
\item The proof of Proposition \ref{Prop_UpperBound_Dens} above is given
 in Appendix \ref{AnnexProofPropUpper}.

\item If $\sigma$ and $b$ is bounded, the classical Aronson's estimates implies that  \eqref{DensUpperBound} holds even without the $|x_{\scriptscriptstyle{0}}|^{4}$ multiplicative term. If $\sigma$ and $b$ are unbounded,  \cite{ArnaudoFeng} provides an adaptation of Aronson's estimates; unfortunately  they first considered time-homogeneous coefficients, and also their result does not imply \eqref{DensUpperBound}.

\item If $\sigma$ and $b$ have polynomial growth and  are time-homogeneous, various estimates are  given in \cite{DeMarcoThesis}. However the behavior is of type $\mathcal{O}(t^{-\frac{3}{2}})$ instead of $\mathcal{O}(t^{-\frac{1}{2}})$ when $t\to 0$.
\end{enumerate}
\end{remark}

Let  $Y_{\kappa}$ be a
random variable distributed according to
$u_{\scriptscriptstyle{0,\kappa}}$.  We are interested in the
following result.

\begin{proposition}\label{Proposition(A)}
Assume that $m \in ]\frac{3}{5},1[$. Let $B$ be a classical Brownian motion independent of $Y_{\kappa}$. Then  there exists a unique
(strong) solution $\overline{Y}=(\overline{Y_t})_{t\in[0,T-\kappa]}$  of
\begin{equation}\label{SDE_Prop(A)}
    \left\{
 \begin{array}{ccl}
   \overline{Y_t}&=&Y_{\kappa}+\int\limits_{0}^t \Phi(\overline{\mathcal{U}}(s,\overline{Y_s}))dB_s,\\
       \overline{\mathcal{U}}(t,\cdot)&=& \mbox{Law  density of}
\ \overline{Y_t},\ \ \forall~t\in [0,T-\kappa],\\
       \overline{\mathcal{U}}(0,\cdot)&=& u_{\scriptscriptstyle{0,\kappa}}\ .
 \end{array}
\right.
\end{equation}
In particular  pathwise uniqueness holds.
\end{proposition}

\begin{corollary}\label{Corollary(A)}

Let $W$ be a classical Brownian motion independent of
$Y_{\kappa}$. Therefore there is a unique (strong) solution
${Y^{\kappa}}=(Y^{\kappa}_t)_{\scriptscriptstyle{t\in[\kappa,T]}}$ of
\begin{equation}\label{SDE_Cor(A)}
    \left\{
 \begin{array}{ccl}
   Y^{\kappa}_t&=&Y_{\kappa}+\int\limits_{\kappa}^t \Phi(\mathcal{U}(s,Y^{\kappa}_s))dW_s,\\
       \mathcal{U}(t,\cdot)&=& \mbox{Law  density of}
\ Y^{\kappa}_t,\ \ \forall~t\in[ \kappa,T],\\
       \mathcal{U}(\kappa,\cdot)&=& u_{\scriptscriptstyle{0,\kappa}}\ .
 \end{array}
\right.
\end{equation}
\end{corollary}

\begin{proof}[Proof of Corollary \ref{Corollary(A)}]

We start with the proof of uniqueness. Let $\kappa >0$. We consider two solutions $Y^{\kappa,1}$ and $Y^{\kappa,2}$ of \eqref{SDE_Cor(A)},  we set $\overline{Y_t^i}=Y^{\kappa,i}_{t+\kappa}$, $\forall t\in[ 0, T-\kappa]$, $i=1,2$ and $B_t=W_{t+k}-W_t$, $\forall t\in[ 0, T-\kappa]$. Clearly  $\overline{Y_t^1}$ and $\overline{Y_t^2}$ solve \eqref{SDE_Prop(A)}. Therefore, using Proposition \ref{Proposition(A)}, we deduce uniqueness for problem \eqref{SDE_Cor(A)}. Existence follows by similar arguments.
\end{proof}

\begin{proof}[Proof of Proposition \ref{Proposition(A)}]

Let $W$ be a classical Brownian motion on some filtered
probability space. Given the function $\overline{\mathcal{U}}$, defined in \eqref{Trans_U}, we
 construct below a unique process $\overline{Y}$  strong
solution of
\begin{equation}\label{Prop(A).(1)}
    \overline{Y_t}=\overline{Y_{\scriptscriptstyle{0}}}+\int\limits_{0}^t \Phi(\overline{\mathcal{U}}(s,\overline{Y_s}))dW_s.
\end{equation}
From
\eqref{Trans_U}, for every $(s,y) \in [0,T-\kappa]\times \mathbb{R}$, we have
\begin{equation*}
    \Phi(\overline{\mathcal{U}}(s,y))=\sqrt{2\bar{a}(s,y)},
\end{equation*}
where
\begin{equation}\label{a_bar}
\bar{a}(s,y)=(s+\kappa)^{\scriptscriptstyle{\alpha(1-m)}}(D+\tilde{k}|y|^2(s+\kappa)^{\scriptscriptstyle{-2\alpha}}).
\end{equation}
In fact,  $\Phi(\overline{\mathcal{U}})$ is continuous, smooth with respect to the space parameter and all the space derivatives  of order greater or equal than one are bounded; in particular $\Phi(\overline{\mathcal{U}})$
is Lipschitz and it has
linear growth. Therefore \eqref{Prop(A).(1)} admits a strong solution.

By Lemma \ref{LemmaI1}  the function $t\mapsto \rho(t,\cdot)$ from $[0,T-\kappa]$ to $\mathcal{M}(\mathbb{R})$, where $\rho(t,\cdot)$ is the law of
$\overline{Y_t}$, is a solution to
\begin{equation}\label{Prop(A).PDE}
    \left\{
 \begin{array}{ccl}
  \partial_t\rho&=& \partial_{xx}^2(\bar{a}\rho),\\
        \rho(0,\cdot)&=&u_{\scriptscriptstyle{0,\kappa}}. \\
\end{array}
\right.
\end{equation}
To conclude it remains to prove that $\overline{\mathcal{U}}(t,y)dy$ is the law of $\overline{Y_t}$, $\forall t\in[0,T-\kappa]$; in particular the law of the r.v. $\overline{Y_t}$ admits a density. For this we will apply Theorem \ref{TheoBRREtendu}  for which we need to check the validity of Hypothesis(B2) when $a=\bar{a}$ and for $z:=z_1-z_2$, where $z_1:=\rho$ and $z_2:=\overline{\mathcal{U}}$. By additivity  this will be of course fulfilled  if we prove it separately for $z:=\rho$ and $z:=\overline{\mathcal{U}}$, which are both solutions to \eqref{Prop(A).PDE}.

Since $\bar{a}$ is non-degenerate, by Remark \ref{RemarkConsq}(1),
we only need to check items (ii) and (iii) of the mentioned Hypothesis(B2). On  one hand,  since $\bar{a}(s,y)=\overline{\mathcal{U}}^{m-1}(s,y)$,  $z:=\overline{\mathcal{U}}$  verifies Hypothesis(B2) because of items (ii) and (iii) of Lemma \ref{Lemma(A1)}. On the other hand, since $\sqrt{{\overline{a}}}$ has linear growth, by
Remark \ref{RemarkConsq}.(4)  $\rho$ fulfills
item (ii) of Hypothesis(B2). Moreover, by Lemma
\ref{Lemma(A5)} below, $\rho$ also verifies
item (iii) of Hypothesis(B2). Finally
Theorem \ref{TheoBRREtendu} implies that
$\overline{\mathcal{U}}\equiv \rho$.

\end{proof}

\begin{lemma}\label{Lemma(A5)}

 Let $\psi:[0,T]\times \mathbb{R}\rightarrow \mathbb{R}_+$, continuous (not necessarily bounded) such that $\psi(t,\cdot)$  is smooth with bounded
derivatives of orders greater or equal than one. We also suppose $\psi$ to be non-degenerate.

We consider a stochastic process $X=(X_t)_{t\in[0,T]}$
strong solution of the SDE
\begin{equation}\label{Lemma(A5).SDE}
    X_t=X_{\scriptscriptstyle{0}}+\int\limits_0^t\psi(s,X_s)dW_s,
\end{equation}
where  $X_{\scriptscriptstyle{0}}$ is a random variable distributed
according to $u_{\scriptscriptstyle{0,\kappa}}$  defined in \eqref{U_kappa} with $m \in ]\frac{3}{5},1[$.

For $t\in]0,T]$  the law of  $X_t$ has a density $\nu(t,\cdot)$ such that  $(\psi^2 \nu)(t,x)$ belongs to $L^2([t_{\scriptscriptstyle{0}},T]\times
\mathbb{R})$, for every $t_{\scriptscriptstyle{0}} >0$.
\end{lemma}

\begin{proof}[Proof of Lemma \ref{Lemma(A5)}]

\({}\)

If $X_0=x_{\scriptscriptstyle{0}}$, where $x_{\scriptscriptstyle{0}}$ is a real number,  then Proposition \ref{Prop_UpperBound_Dens}  implies that,  for every $t\in ]0,T]$, the law of $X_t$ admits a density $p_t(x_{\scriptscriptstyle{0}},\cdot)$. Consequently,  if  the law of $X_0$ is $u_{0,\kappa}(x)dx$, for every $t\in]0,T]$, the law of $X_t$  has a density  given by
\begin{equation*}
    \nu(t,x)=\int\limits_{\mathbb{R}}u_{\scriptscriptstyle{0,\kappa}}(x_{\scriptscriptstyle{0}})p_t(x_{\scriptscriptstyle{0}},x)dx_{\scriptscriptstyle{0}}.
\end{equation*}
By \eqref{DensUpperBound} in Proposition \ref{Prop_UpperBound_Dens}  it follows
\begin{equation}\label{DensUpperBoundBis}
    \sup_{(t,x)\in[t_{\scriptscriptstyle{0}},T]\times
\mathbb{R}}p_t(x_{\scriptscriptstyle{0}},x)\leq K_0(1+|x_{\scriptscriptstyle{0}}|^4),\quad \mbox{ where} \quad K_0=\frac{\displaystyle{K}}{\displaystyle{\sqrt{t_{\scriptscriptstyle{0}}}}}.
\end{equation}
Using \eqref{DensUpperBoundBis}  we get
\begin{equation}\label{Sup_Nu}
    K_1:=\sup_{(t,x)\in[t_{\scriptscriptstyle{0}},T]\times \mathbb{R}}|\nu(t,x)|\leq K_0
    \int\limits_{\mathbb{R}}(1+|x_{\scriptscriptstyle{0}}|^4)\mathcal{U}(\kappa,x_{\scriptscriptstyle{0}})dx_{\scriptscriptstyle{0}}< \infty;
\end{equation}
the latter inequality is valid because of  \eqref{L(A1).(ii).1} in Lemma
\ref{Lemma(A1)}. In the sequel of the proof,  the constants $K_2,K_3,K_4$  will only depend on $t_{\scriptscriptstyle{0}}$, $T$ and $\psi$. Furthermore
\begin{equation*}
\int\limits_{[t_{\scriptscriptstyle{0}},T]\times
\mathbb{R}}\left((\psi^{\scriptscriptstyle{2}}\nu)(t,x)\right)^{\scriptscriptstyle{2}}dtdx\leq
\sup_{(t,x)\in[t_{\scriptscriptstyle{0}},T]\times \mathbb{R}}|\nu(t,x)|\mathds{E}\left[\int\limits_0^T \psi^{\scriptscriptstyle{4}}(t,X_t)dt\right].
\end{equation*}
Since $\psi$ has linear growth, this expression is bounded by
\begin{equation}
    K_1 K_2 \left(1+\int\limits_0^T\mathds{E}\left[ \sup_{t\in[0,T]}|X_t|^{\scriptscriptstyle{4}}\right]dt\right).\label{Lemma(A5).(1)}
\end{equation}
\eqref{Lemma(A5).(1)} follows because of \eqref{Sup_Nu}. Besides, by Burkholder-Davis-Gundy and Jensen's inequalities, taking into account  the linear growth of $\psi$,  it follows that
\begin{equation*}
    \mathds{E}\left[\sup_{t\in[0,T]}|X_t|^{\scriptscriptstyle{4}}\right]\leq
    K_3\left(\mathds{E}\left[|X_{\scriptscriptstyle{0}}|^4\right]+
    \int\limits_0^T\mathds{E}\left[\sup_{s\in[0,T]}|X_s|^{\scriptscriptstyle{4}}\right]ds+T\right).
\end{equation*}
Then, by Gronwall's lemma, there is another constant $K_4$  such that
\begin{equation}\label{Lemma(A5).(2)}
   \mathds{E}\left[\sup_{t\in[0,T]}|X_t|^{\scriptscriptstyle{4}}\right]\leq K_4\left(1+\int\limits_{\mathbb{R}}|x_{\scriptscriptstyle{0}}|^4\mathcal{U}(\kappa,x_{\scriptscriptstyle{0}})dx_{\scriptscriptstyle{0}}\right).
\end{equation}
Finally  \eqref{Lemma(A5).(1)}, \eqref{Lemma(A5).(2)} and \eqref{Sup_Nu}  allow  us to conclude the proof.
\end{proof}

We are now ready to provide the probabilistic representation related
to function $\mathcal{U}$ which in fact is only  a solution in law of \eqref{NLSDE}.

\begin{definition}
We say that \eqref{NLSDE} admits a weak (in law) solution if there is a probability space $(\Omega, \mathcal{F},\mathds{P})$, a Brownian motion $(W_t)_{t\geq 0}$ and a process $(Y_t)_{t\geq 0}$ such that   the system \eqref{NLSDE} holds. \eqref{NLSDE}  admits uniqueness in law if,  given $(W^1,Y^1)$, $(W^2,Y^2)$ solving \eqref{NLSDE} on some related probability space, it follows that $Y^1$ and $Y^2$ have the same law.
\end{definition}

\begin{theorem}\label{Theorem(B)}
Assume that $m\in ]\frac{3}{5},1[$. Then  there is a unique weak
solution (in law)  $Y$ of  problem \eqref{NLSDE}.
\end{theorem}
\begin{remark}
Indeed the assumption on $m\in]\frac{3}{5},1[$ is only required for the application of Theorem \ref{TheoBRREtendu}. The arguments following  the present proof only use $m>\frac{1}{3}$.
\end{remark}
\begin{proof}[Proof of Theorem \ref{Theorem(B)}]

First we start with the existence of a weak solution for
\eqref{NLSDE}. Let $\mathcal{U}$ be again the (Barenblatt's) solution of
\eqref{PME}.  We consider the solution $(Y_t^{\kappa})_{t\in[\kappa,T]}$ provided by Corollary \ref{Corollary(A)} extended to $[0, \kappa]$, setting $Y_t^{\kappa}=Y_{\kappa}$, $ t\in[0,\kappa]$. We prove that the laws of  processes $Y^{\kappa}$ are tight. For this  we implement
the classical Kolmogorov's criterion, see \cite[Section 2.4, Problem 4.11 ]{KarShreBook}. We will show  the existence of $p>2$ such that
\begin{equation}\label{Theo(B).(1)}
    \mathds{E}\left[\left|Y_t^{\kappa}-Y_s^{\kappa}\right|^{\scriptscriptstyle{p}}\right]\leq
    \mathcal{C}_p|t-s|^{\scriptscriptstyle{\frac{p}{2}}},\ \ \ \forall s,t \in [0,T],
\end{equation}
where $\mathcal{C}_p$ will stand for a constant (not always the same),
  depending on $p$ and $T$ but not on $\kappa$. Let $s,t \in ]0,T]$. Let $p>2$. By Burkholder-Davis-Gundy inequality  we obtain
\begin{equation*}
    \mathds{E}\left[\left|Y_t^{\kappa}-Y_s^{\kappa}\right|^{\scriptscriptstyle{p}}\right]\leq
    \mathcal{C}_p
    \mathds{E}\left[\left|\int\limits_s^t\Phi^2(\mathcal{U}(r,Y_r^{\kappa}))dr\right|^{\scriptscriptstyle{\frac{p}{2}}}\right].
\end{equation*}
Then, using Jensen's inequality and the fact that $\mathcal{U}(r,\cdot)$ is the law density of $Y_r^{\kappa}$, $r\geq \kappa$, we get
\begin{equation}\label{Theo(B).(2)}
    \mathds{E}\left[\left|Y_t^{\kappa}-Y_s^{\kappa}\right|^{\scriptscriptstyle{p}}\right]\leq
   \mathcal{C}_p|t-s|^{\scriptscriptstyle{\frac{p}{2}-1}}
    \int\limits_s^t dr\int\limits_{\mathbb{R}}\Phi^p(\mathcal{U}(r,y))\mathcal{U}(r,y)dy.
\end{equation}
We have
\begin{equation*}
\int\limits_s^t
dr\int\limits_{\mathbb{R}}\Phi^p(\mathcal{U}(r,y))\mathcal{U}(r,y)dy=\int\limits_s^t
dr\int\limits_{\mathbb{R}}dy\left(\mathcal{U}(r,y)\right)^{\scriptscriptstyle{\frac{p(m-1)}{2}+1}},\nonumber\\
\end{equation*}
and,  by Lemma \ref{Lemma(A1)} (i), the result follows.

Consequently  there is a subsequence $Y^n:=Y^{\kappa_n}$  converging
in law (as $C([0,T])-$valued random elements) to some process $Y$. Let $P^n$ be the corresponding laws on the canonical space $\Omega=C([0,T])$ equipped with the Borel $\sigma$-field. $Y$ will denote the canonical process $Y_t(\omega)=\omega(t)$. Let $P$ be  the weak limit of $(P^n)$.

1) We first observe that the marginal laws of $Y$ under $P^n$ converge to the marginal law of $Y$ under $P$. Let $t\in]0,T]$. If the sequence $(\kappa_n)$ is lower than $t$, then the law of $Y_t$ under $P^n$ equals the constant law $\mathcal{U}(t,x)dx$. Consequently, for every $t\in]0,T]$, the law of $Y_t$ under $P$ is $\mathcal{U}(t,x)dx$.

2) We now prove   that $Y$ is a (weak) solution of \eqref{NLSDE},  under $P$. By similar arguments as for the classical stochastic differential equations, see \cite[Chapter 6]{StrVarBook}, it is enough to prove that $Y$ (under $P$) fulfills the martingale problem i.e., for every $f\in C^2_b(\mathbb{R})$, the process
\begin{equation*}
\mbox{\textbf{(MP)}} \ \ \ \
f(Y_t)-f(0)-\frac{\displaystyle{1}}{\displaystyle{2}}\int\limits_0^tf''(Y_s)\Phi^2(\mathcal{U}(s,Y_s))ds,
\end{equation*}
is an $(\mathcal{F}_s)$-martingale, where $(\mathcal{F}_s)$ is the canonical filtration associated with $Y$. $C^2_b(\R)$ stands for the set $\{f\in C^2(\R)\vert f,f',f'' \mbox{ bounded }\}$.  Let $\mathds{E}$ (resp. $\mathds{E}^n$) be the expectation operator with respect to $P$ (resp. $P^n$). Let $s,t \in [0,T]$  with $s<t$ and $R=R(Y_r,r\leq s)$ be an
$\mathcal{F}_s-$measurable, bounded and continuous (on $C([0,T])$) random variable. In order to show the martingale
property \textbf{(MP)} of $Y$, we have to prove that
\begin{equation}\label{MP_Bis}
\mathds{E}\left[\left(f(Y_t)-f(Y_s)-\frac{\displaystyle{1}}{\displaystyle{2}}\int\limits_s^tf''(Y_r)\Phi^2(\mathcal{U}(r,Y_r))dr\right)R\right]=0, \ \ \ f\in C^2_b(\mathbb{R}).
\end{equation}
We first consider the case when $s>0$. There is $n\geq n_{\scriptscriptstyle{0}}$, such that $\kappa_n<s$. Let $f\in C^2_b(\mathbb{R})$;
 since $(Y_s)_{s\geq \kappa_n}$,   under $P^n$,  are still martingales  we have
\begin{equation}\label{Theo(B).(5)}
\mathds{E}^n\left[\left(f(Y_t)-f(Y_s)-\frac{\displaystyle{1}}{\displaystyle{2}}\int\limits_s^tf''(Y_r)\Phi^2(\mathcal{U}(r,Y_r))dr\right)R\right]=0.
\end{equation}
We are able to prove that \eqref{MP_Bis} follows from \eqref{Theo(B).(5)}. Let $\varepsilon>0$ and $N>0$ such that
\begin{equation}\label{D10}
\int\limits_s^tdr\int\limits_{\{|y|>\frac{N}{C}-1\}}\mathcal{U}^m(r,y)dy\leq\varepsilon,
\end{equation}
where  $C$ is the linear growth constant of $\Phi^2\circ\mathcal{U}$ in the sense of Definition \ref{LinGrowth}. In order to conclude,  passing to the limit in \eqref{Theo(B).(5)}, we will only have to show that
\begin{equation}\label{D11}
    \lim_{n\to+\infty}\mathds{E}^n\left[F(Y)\right]-\mathds{E}\left[F(Y)\right]=0,
\end{equation}
where
   $ F(\ell)=\int\limits_s^tdr\Phi^2(\mathcal{U}(r,\ell(r)) f''(\ell(r))  R(\ell(\xi),\xi\leq s)$,
$F:C([0,T])\rightarrow\R$ being continuous but not bounded. The left-hand side of \eqref{D11} equals
\begin{align}
   & \mathds{E}^n\left[F(Y)-F^N(Y)\right]+\mathds{E}^n\left[F^N(Y)\right]-\mathds{E}\left[F^N(Y)\right]+\mathds{E}\left[F^N(Y)-F(Y)\right]\nonumber\\
   \label{D15}\\
   &:=\mathcal{E}_1(n,N)+\mathcal{E}_2(n,N)+\mathcal{E}_3(n,N),\nonumber
\end{align}
where
\begin{equation*}
    F^N(\ell)=\int\limits_s^tdr\left(\Phi^2(\mathcal{U}(r,\ell(r))\wedge N\right)f''(\ell(r))R(\ell(\xi),\xi\leq s).
\end{equation*}
Since $\kappa_n<s$, for $N$ large enough,  we get
\begin{align}
|\mathcal{E}_1(n,N)|&\leq \|R\|_{\infty}\|f''\|_{\infty}\int\limits_s^tdr\int\limits_{\{\Phi^2(\mathcal{U}(r,y))\geq N\}}\left(\Phi^2(\mathcal{U}(r,y)-N\right)\mathcal{U}(r,y)dy \nonumber\\
&\leq  \|R\|_{\infty}\|f''\|_{\infty}\int\limits_s^tdr\int\limits_{\{|y|>\frac{N}{C}-1\}}\mathcal{U}^m(r,y)dy\leq \varepsilon\|R\|_{\infty}\|f''\|_{\infty},\label{D14}
\end{align}
taking into account \eqref{D10} and the second item of Lemma \ref{Lemma(A1)}. For fixed $N$,  chosen in \eqref{D10},   we have
$\lim_{n\to+\infty}\mathcal{E}_2(n,N)=0$,
since $F^N$ is bounded and continuous. Again, since the law density under $P$ of $Y_t$, $t\geq s$ , is $\mathcal{U}(t,\cdot)$, similarly as for \eqref{D14}, we obtain $|\mathcal{E}_3(n,N)|\leq \varepsilon\|R\|_{\infty}\|f''\|_{\infty}$. Finally, coming back to \eqref{D15}, it follows
\begin{equation*}
    \limsup_{n\to+\infty} \left|\mathds{E}^n\left[F(Y)\right]-\mathds{E}\left[F(Y)\right]\right|\leq 2\varepsilon\|R\|_{\infty}\|f''\|_{\infty};
\end{equation*}
since $\varepsilon>0$ is arbitrary, \eqref{D11} is established. So  \eqref{MP_Bis} is verified for $s>0$.

3) Now, we consider the case when $s=0$. We first prove that
\begin{equation}\label{(E1).1}
    \mathds{E}\left[ \int\limits_0^T \Phi^2(\mathcal{U}(r,Y_r))dr\right]<+\infty.
\end{equation}
By item 1) of this proof, the law of $Y_r$, $r>0$ admits $\mathcal{U}(r,\cdot)$ as density. Consequently, the left-hand side of \eqref{(E1).1} gives
\begin{equation*}
    \int\limits_{]0,T]}dr\int\limits_{\mathbb{R}} \Phi^2(\mathcal{U}(r,y))\mathcal{U}(r,y)dy=\int\limits_0^Tdr\int\limits_{\mathbb{R}}\mathcal{U}^m(r,y)dy,
\end{equation*}
which is finite  by the second item of Lemma \ref{Lemma(A1)}. Coming back to \eqref{MP_Bis}, we can now let $s$ go to zero. Since  $Y$ is continuous and $f$ is bounded,  we clearly have $\lim_{s\to 0}\mathds{E}\left[f(Y_s)R\right]=\mathds{E}\left[f(Y_0)R\right]$. Moreover
\begin{equation*}
    \lim_{s\to 0}\mathds{E}\left[\left(\int\limits_s^tf''(Y_r)\Phi^2(\mathcal{U}(r,Y_r))dr\right)R\right]
    =\mathds{E}\left[\left(\int\limits_0^tf''(Y_r)\Phi^2(\mathcal{U}(r,Y_r))dr\right)R\right],
\end{equation*}
using Lebesgue's dominated convergence theorem and \eqref{(E1).1}. Consequently  we obtain
\begin{equation}\label{(E1).2}
\mathds{E}\left[\left(f(Y_t)-f(Y_0)-\frac{\displaystyle{1}}{\displaystyle{2}}\int\limits_0^tf''(Y_r)\Phi^2(\mathcal{U}(r,Y_r))dr\right)R\right]=0.
\end{equation}
It remains to show that $Y_0=0$ a.s.  This follows because $Y_t \rightarrow Y_0$  a.s., and also in law (to $\delta_0$),  by the third item of Proposition \ref{Prop_U_pp}. Finally  we have shown that the limiting process $Y$ verifies \textbf{(MP)}, which proves  the existence of solutions to \eqref{NLSDE}.

4) We now prove  uniqueness. Since $\mathcal{U}$ is fixed,  only  uniqueness for the first line of equation \eqref{NLSDE} has to be  established. Let $(Y_t^{i})_{t\in[0,T]}$, $i=1,2$, be two solutions. In order to show that the laws of $Y^1$ and $Y^2$ are identical, according to  \cite[Lemma 2.5]{Jacod_SFlour},  we will verify that their finite marginal distributions are the same. For this  we consider $0=t_{\scriptscriptstyle{0}}<t_{\scriptscriptstyle{1}}<\ldots<t_{\scriptscriptstyle{N}}=T$. Let $0<\kappa<t_{\scriptscriptstyle{1}}$. Obviously  we have $Y^i_{\scriptscriptstyle{t_0}}=0$ a.s., in the corresponding probability space, $\forall i  \in \{1,2\}$.  Both restrictions $Y^1|_{[\kappa,T]}$ and $Y^2|_{[\kappa,T]}$ verify \eqref{SDE_Cor(A)}. Since that equation admits pathwise uniqueness, it also admits uniqueness in law by Yamada-Watanabe theorem. Consequently $Y^1|_{[\kappa,T]}$ and $Y^2|_{[\kappa,T]}$ have the same law  and in conclusion the law  of $(Y_{t_1}^1,\ldots,Y_{t_N}^1)$ coincides with the law of $(Y_{t_1}^2,\ldots,Y_{t_N}^2)$, thus, the law of $(Y_{t_0}^1,\ldots,Y_{t_N}^1)$ coincides with the law of $(Y_{t_0}^2,\ldots,Y_{t_N}^2)$.
\end{proof}


\section{Numerical experiments}\label{NumExp}
In order to avoid singularity problems due to the initial condition
being a Dirac delta function, we will consider a time translation of
$\mathcal{U}$, denoted $v$, and defined by
\[v(t,\cdot)=\mathcal{U}(t+1,\cdot),\ \ \forall t \in [0,T].\]
 $v$ still solves equation \eqref{PME}, for $m\in]0,1[$,
but with now  a smooth initial data given by
\begin{equation}\label{V_0}
    v_{\scriptscriptstyle{0}}(x)=\mathcal{U}(1,x),\ \ \ \forall x\in \mathbb{R}.
\end{equation}
Indeed, we have the following formula
\begin{equation}\label{FormExpTrans}
    v(t,x)=(t+1)^{-\alpha}{\left(D+\tilde{k}|x|^2(t+1)^{-2\alpha}\right)}^{-\frac{1}{1-m}},
\end{equation}
where $\alpha$, $\tilde{k}$ and $D$ are still given by
\eqref{Constants}.

We now wish  to compare the exact solution of  problem \eqref{PME}
to a numerical probabilistic solution. In fact, in order to perform
such approximated solutions, we use the algorithm described in
Sections 4  of \cite{BCR1} (implemented in Matlab).
 We focus on the case $m=\frac{\displaystyle{1}}{\displaystyle{2}}$. \\

\begin{figure}[htbp]
\begin{center}
\includegraphics[width=\linewidth,height=0.33\textheight]{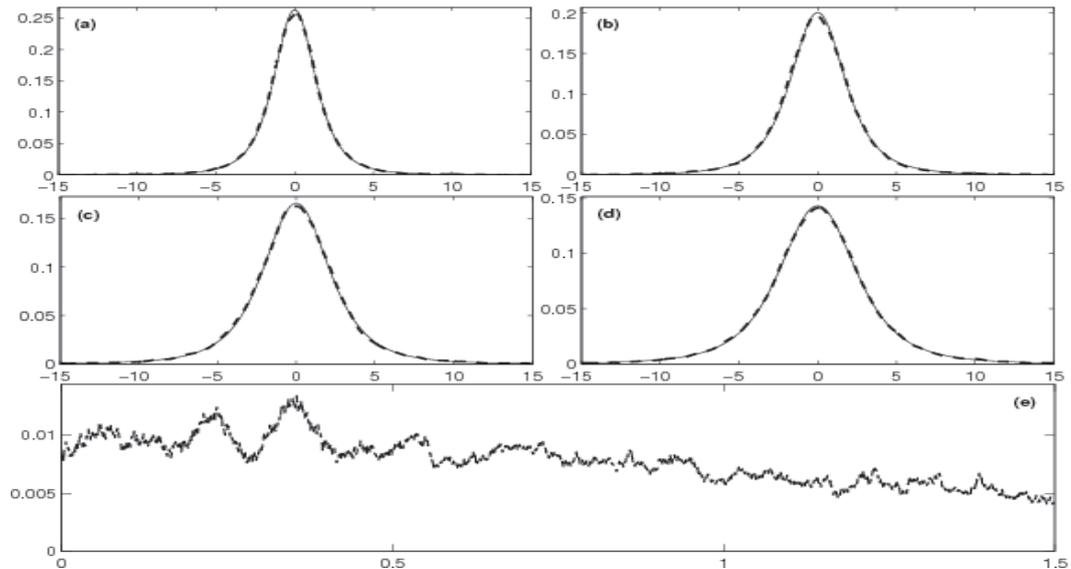}
\caption{Numerical (dashed line) and exact solutions (solid line) values
at t=0 (a),
   t=0.5 (b), t=1 (c) and t=1.5 (d). The evolution of the $L^2$ error
    over the time interval $[0,1.5]$
   (e).}
\label{fig:PMEm_0.5_PROB_Sol_Err}
\end{center}
\end{figure}

\textbf{Simulation experiments:} we compute the
numerical solution over the time-space grid
 $[0,1.5]\times[-15,15]$. We use  $n=50000$ particles and a time step
    $\Delta t=2\times 10^{-4}$. Figures  \ref{fig:PMEm_0.5_PROB_Sol_Err}.(a)-(b)-(c)-(d), display the
exact and the numerical solutions at times $t=0$,
$t=0.5$, $t=1$ and $t=T=1.5$, respectively. The exact solution for
the fast diffusion equation \eqref{PME}, given in
\eqref{FormExpTrans}, is depicted by solid lines. Besides, Figure \ref{fig:PMEm_0.5_PROB_Sol_Err}.(e) describes the
time evolution of  the discrete $L^2$ error on the
time interval $[0,1.5]$.

\section{Appendix}\label{Annex}

\subsection{Proof of Proposition \ref{Prop_UpperBound_Dens}}\label{AnnexProofPropUpper}

We start with some notations for the Malliavin calculus. The set $\mathds{D}^{\infty}$  represents the classical Sobolev-Malliavin
space of smooth test random variables. $\mathds{D}^{1,2}$ is defined in the lines after  \cite[Lemma 1.2.2]{LivreNual}
 and $\mathds{L}^{1,2}$ is introduced in   \cite[Definition 1.3.2]{LivreNual}. See also \cite{LivreMalliavin} for a complete monograph on Malliavin calculus. We state a preliminary result.
\begin{proposition}\label{Prop_E(N{-p})}

Let $N$ be a non-negative random variable. Suppose, for
every $p\geq 1$, the  existence of
 constants $C(p)$ and $\epsilon_{\scriptscriptstyle{0}}(p)$ such that
\begin{equation}\label{Cdt}
    {P}(N\leq \epsilon)\leq C(p).\epsilon^{p+1},\ \ \
    \forall \epsilon\in]0,\epsilon_{\scriptscriptstyle{0}}(p)].
\end{equation}
Then, for every $p\geq 1$,
\begin{equation}\label{E(N{-p})}
\mathds{E}(N^{-p})\leq
\epsilon_{\scriptscriptstyle{0}}(p).C(p+1)+\epsilon_{\scriptscriptstyle{0}}(p)^{-p}{P}(N>
\epsilon_{\scriptscriptstyle{0}}(p)).
\end{equation}
\end{proposition}

\begin{proof}[Proof of Proposition \ref{Prop_E(N{-p})}]
Let $p\geq 1$ and $\epsilon_{\scriptscriptstyle{0}}(p)>0$. Setting  $F(x)={P}(N\leq x)$,  $ x\in\mathbb{R}_+$,
we have
\begin{equation}\label{Prop_E(N{-p}).(1)}
\mathds{E}(N^{-p})=I_1+I_2,
\end{equation}
where
\begin{equation*}
    I_1=\int\limits_0^{\epsilon_{\scriptscriptstyle{0}}(p)}x^{-p}dF(x) \mbox{ and } I_2=\int\limits_{\epsilon_{\scriptscriptstyle{0}}(p)}^{+\infty}x^{-p}dF(x).
\end{equation*}
\eqref{Cdt} implies that $I_1$ and $I_2$ are well-defined. Indeed, on one hand, applying integration by parts on $I_1$, we get
\begin{equation*}
I_1={\left[x^{-p}F(x)\right]}^{\epsilon_{\scriptscriptstyle{0}}(p)}_0+p\int\limits_0^{\epsilon_{\scriptscriptstyle{0}}(p)}x^{-p-1}F(x)dx;
\end{equation*}
moreover,  there is a constant $C(p)$ such that
\begin{equation}\label{Prop_E(N{-p}).(2)}
I_1\leq (p+1)\epsilon_{\scriptscriptstyle{0}}(p)C(p).
\end{equation}
On the other hand, again \eqref{Cdt} says that
\begin{equation}\label{Prop_E(N{-p}).(3)}
  I_2\leq  {\epsilon_{\scriptscriptstyle{0}}(p)}^{-p}(1-F(\epsilon_{\scriptscriptstyle{0}}(p))).
\end{equation}
Consequently, using
\eqref{Prop_E(N{-p}).(2)} and \eqref{Prop_E(N{-p}).(3)} and coming back to \eqref{Prop_E(N{-p}).(1)},
\eqref{E(N{-p})} is established.
\end{proof}

\begin{proof}[Proof of Proposition \ref{Prop_UpperBound_Dens}]

In this proof $\sigma'$ (resp. $b'$) stands for $\partial_x \sigma$ (resp. $\partial_x b$). Let $Y =(Y_t^{x_{\scriptscriptstyle{0}}})_{t\in[0,T]}$, be  the
solution of \eqref{SDE_Upper_Bound}.  According to  \cite[Theorem 2.2.2]{LivreNual}  we have  $Y_s
\in \mathds{D}^{\infty}$, $\forall s \in[0,T]$. Let $s>0$. Since $\sigma$ is non-degenerate, by
 \cite[Theorem 2.3.1]{LivreNual}, the law of $Y_s$ admits a density that we denote
by $p_s(x_{\scriptscriptstyle{0}},\cdot)$.

The second step consists in a re-scaling, transforming the time $s$ into a noise multiplicative  parameter $\lambda$; we set $\lambda=\sqrt{s}$. Indeed, $(Y_t)$ is distributed as $(Y^{\lambda}_{\frac{t}{\lambda^2}})$, where
\begin{equation*}
    Y_t^{\lambda}=x_{\scriptscriptstyle{0}}+\lambda\int\limits_0^t\sigma(r\lambda^2,Y_r^{\lambda})dW_r+\lambda^{\scriptscriptstyle{2}}\int\limits_0^tb(r\lambda^2,Y_r^{\lambda})dr.
\end{equation*}
In particular, $Y_s \sim Y_1^{\lambda}$. By previous arguments, for every $t>0$, $Y_t^{\lambda}
\in \mathds{D}^{\infty}$ and its law  admits a density denoted by $p_t^{\lambda}(x_{\scriptscriptstyle{0}},\cdot)$. Our aim consists in showing the existence of a constant $K$  such that
\begin{equation}\label{(A1)}
    p_1^{\lambda}(x_{\scriptscriptstyle{0}},y)\leq \frac{\displaystyle{K}}{\displaystyle{\lambda}}(1+|x_{\scriptscriptstyle{0}}|^4),\ \forall y\in\mathbb{R}, \lambda \in ]0,\sqrt{T}],
\end{equation}
where  $K$ is a constant which does not depend on $x_{\scriptscriptstyle{0}}$ and $\lambda$. In fact, we will prove that, for every $\lambda \in ]0,\sqrt{T}]$,
\begin{equation}\label{A11}
     \sup_{y\in\mathbb{R}, t\in ]0,1]}p_t^{\lambda}(x_{\scriptscriptstyle{0}},y)\leq \frac{\displaystyle{K}}{\displaystyle{\lambda}}(1+|x_{\scriptscriptstyle{0}}|^4).
\end{equation}
We set $ Z_t^{\lambda}=\frac{\displaystyle{Y_t^{\lambda}-x_{\scriptscriptstyle{0}}}}{\displaystyle{\lambda}}, \ \ t \in[0,1]$,
so that  the density $q_t^{\lambda}$ of $Z_t^{\lambda}$ fulfills $ q_t^{\lambda}(z)=\lambda p_t^{\lambda}(x_{\scriptscriptstyle{0}}, \lambda z+ x_{\scriptscriptstyle{0}}), \  (t,z)\in [0,1]\times \mathbb{R}$.
In fact, we will have attained \eqref{A11}, if we show
\begin{equation}\label{(A2)}
    \sup_{z\in\mathbb{R}, \lambda \in ]0,\sqrt{T}]}q_t^{\lambda}(z)\leq K(1+|x_{\scriptscriptstyle{0}}|^4), \ t\in]0,1].
\end{equation}
We express the equation fulfilled by $Z$; it yields
\begin{equation}\label{(A3)}
    Z_t^{\lambda}=\int\limits_0^t\sigma^{\lambda}(r,Z_r^{\lambda})dW_r+\int\limits_0^tb^{\lambda}(r,Z_r^{\lambda})dr,
\end{equation}
where,  for every $(r,z)\in [0,1]\times \mathbb{R}$, we set
\begin{equation*}
    \sigma^{\lambda}(r,z)=\sigma(r\lambda^{\scriptscriptstyle{2}},\lambda z+x_{\scriptscriptstyle{0}}), \ \mbox{ and } b^{\lambda}(r,z)=\lambda b(r\lambda^{\scriptscriptstyle{2}},\lambda z+x_{\scriptscriptstyle{0}}).
\end{equation*}
At this stage we state the following lemma.

\begin{lemma}\label{Lemma(A3)}

For $\lambda \in ]0,1]$, we shorten by $Z:=(Z_t^{\lambda})_{t\in[0,T]}$,   the
solution of \eqref{(A3)}. For every $\gamma \geq
1$, we have
\begin{equation}\label{Sup_Z_t}
   \sup_{\lambda \in ]0,\sqrt{T}]} \mathds{E}\left[\sup_{t\in[0,1]}|Z_t|^{\gamma}\right]\leq
    C(1+|x_{\scriptscriptstyle{0}}|^{\gamma}),
\end{equation}
where  $C$ is a constant depending on $\|\sigma'\|_{\infty}$,
$\|b'\|_{\infty}$ and $T$, but  not  on
$x_{\scriptscriptstyle{0}}$.
\end{lemma}

\begin{remark}\label{RemarkNotation}
\({}\)
\begin{enumerate}
  \item For simplicity,  in the whole proof of Proposition \ref{Prop_UpperBound_Dens}, we will set $T=1$.

  \item Since there is no more ambiguity,  we will use again the letter $s$ in the considered integrals.
\end{enumerate}
\end{remark}

\begin{proof}[Proof of Lemma \ref{Lemma(A3)}]

Let $\lambda \in ]0,1]$ and $\gamma \geq 1$.  In the proof $C_1$ is a constant depending on $T$, and $C_2$, $C_3$ depend on $T$, $\|\sigma'\|_{\infty}$ and
$\|b'\|_{\infty}$. Using Burkholder-Davis-Gundy and Jensen's inequalities, we get
\begin{equation*}
     \mathds{E}\left[\sup_{\rho\in[0,t]}|Z_{\rho}|^{\gamma}\right]\leq C_1\left( \int_0^t\mathds{E}\left[|\sigma(s\lambda^{\scriptscriptstyle{2}},\lambda Z_s+ x_{\scriptscriptstyle{0}})|^{\gamma}\right]ds+ \int_0^t\mathds{E}\left[|\lambda b(s\lambda^{\scriptscriptstyle{2}},\lambda Z_s+ x_{\scriptscriptstyle{0}})|^{\gamma}\right]ds\right).
\end{equation*}
Since $\sigma'$ and $b'$ are bounded,  $\sigma$ and
$b$ have linear growth. Therefore, previous expression is bounded by
\begin{equation*}
      C_2(1+\lambda^{\gamma})\left(
     1+|x_{\scriptscriptstyle{0}}|^{\gamma}+\lambda^{\gamma}
     \int_0^t\mathds{E}\left[\sup_{\rho\in[0,s]}|Z_{\rho}|^{\gamma}\right]ds\right).
\end{equation*}
Since $\lambda \in]0,1]$  we obtain
\begin{equation*}
     \mathds{E}\left[\sup_{\rho\in[0,t]}|Z_{\rho}|^{\gamma}\right]\leq C_3\left(
     1+|x_{\scriptscriptstyle{0}}|^{\gamma}+
     \int_0^t\mathds{E}\left[\sup_{\rho\in[0,s]}|Z_{\rho}|^{\gamma}\right]ds\right),
\end{equation*}
Consequently, using Gronwall's lemma, the
result follows.
\end{proof}

Now, in order to perform \eqref{(A2)}, we make use of Malliavin calculus for deriving  expression \eqref{(A3)}. Omitting $\lambda$ in the notation $Z_t^{\lambda}$, we get
\begin{equation*}
    D_rZ_t=\sigma(r\lambda^{\scriptscriptstyle{2}},\lambda Z_r+x_{\scriptscriptstyle{0}} )\mathds{1}_{[r,1]}(t)+\lambda\int\limits_r^t\sigma'(s\lambda^{\scriptscriptstyle{2}},\lambda Z_s+x_{\scriptscriptstyle{0}})D_rZ_sdW_s+\lambda^{\scriptscriptstyle{2}}\int\limits_r^tb'(s\lambda^{\scriptscriptstyle{2}},\lambda Z_s+x_{\scriptscriptstyle{0}})D_rZ_sds.
\end{equation*}
Consequently
\begin{equation*}
    D_rZ_t=\sigma(r\lambda^{\scriptscriptstyle{2}},\lambda Z_r+x_{\scriptscriptstyle{0}} )\mathcal{E}\left(\lambda\int\limits_r^t\sigma'(s\lambda^{\scriptscriptstyle{2}},\lambda Z_s+x_{\scriptscriptstyle{0}})dW_s+\lambda^{\scriptscriptstyle{2}}\int\limits_r^tb'(s\lambda^{\scriptscriptstyle{2}},\lambda Z_s+x_{\scriptscriptstyle{0}})ds\right), \ r<t,
\end{equation*}
where  $\mathcal{E}(S)$ denotes the Dol\'eans exponential of the
continuous semi-martingale
\begin{equation*}
    S_t=\lambda\int\limits_r^t\sigma'(s\lambda^{\scriptscriptstyle{2}},\lambda Z_s+x_{\scriptscriptstyle{0}})dW_s+\lambda^{\scriptscriptstyle{2}}\int\limits_r^tb'(s\lambda^{\scriptscriptstyle{2}},\lambda Z_s+x_{\scriptscriptstyle{0}})ds, \ t\in[0,1].
\end{equation*}
We recall that, $\ \mathcal{E}_t(S)=\exp(S_t-\frac{1}{2}[S]_t)$. Consequently, for fixed $t\in]0,1]$, we have
\begin{equation}\label{Prop_UpperBound_Dens.(1)}
    <DZ_t,DZ_t>=\int\limits_0^t\sigma^2(r\lambda^{\scriptscriptstyle{2}},\lambda Z_r+x_{\scriptscriptstyle{0}})\mathcal{E}^2(\lambda;r,t)dr,
\end{equation}
where
\begin{equation}\label{Prop_UpperBound_Dens.(1bis)}
    \mathcal{E}(\lambda;r,t)=\exp\left(\lambda\int\limits_r^t\sigma'(s\lambda^{\scriptscriptstyle{2}},\lambda Z_s+x_{\scriptscriptstyle{0}})dW_s+\lambda^{\scriptscriptstyle{2}}\int\limits_r^t(b'-\frac{(\sigma')^2}{2})(s\lambda^{\scriptscriptstyle{2}},\lambda Z_s+x_{\scriptscriptstyle{0}})ds\right).
\end{equation}
We set, for every $s\leq t$, $ G(\lambda,s)=\frac{\displaystyle{D_sZ_t}}{\displaystyle{<DZ_t,DZ_t>}}$.
 In view of the application of
\cite[Proposition 2.1.1]{LivreNual}, which implies the useful expression \eqref{Prop_UpperBound_Dens.(3bisbis)} for the density of $Z_t$, we will need to show that $G(\lambda,\cdot)$ belongs to the
domain of the divergence operator $\delta$, denoted by  $Dom \
\delta$. It will be the case if $G(\lambda,\cdot) \in \mathds{L}^{1,2}(H)$  with
$H=L^2([0,T])$. In fact, by the lines after
\cite[Definition 1.3.2]{LivreNual}, we know that $\mathds{L}^{1,2}\subset Dom \
\delta$. Since $Z_t \in \mathds{D}^{\scriptscriptstyle{\infty}}$, we can deduce that
$\frac{\displaystyle{1}}{\displaystyle{<DZ_t,DZ_t>}}$ belongs to
$\mathds{D}^{\scriptscriptstyle{\infty}}$, provided that we prove
\begin{equation}\label{(C1)}
    \frac{\displaystyle{1}}{\displaystyle{<DZ_t,DZ_t>}}\in
L^p(\Omega), \ \ \forall p\geq 1,
\end{equation}
see  \cite[Lemma 2.1.6]{LivreNual}. Since $\mathds{D}^{\scriptscriptstyle{\infty}}$ is an algebra,  $G(\lambda,s) \in \mathds{D}^{\scriptscriptstyle{\infty}}$,  for $s\in ]0,T]$,  and so $G(\lambda, s)\in \mathds{D}^{\scriptscriptstyle{1,2}}$. \eqref{(C1)}  will be the object of Proposition \ref{Prop(C)}.  According to   \cite[Definition 1.3.2]{LivreNual}, to affirm that $G(\lambda,\cdot)$ belongs to $\mathds{L}^{1,2}$ it remains
to show the existence of a measurable version of $D_{s_1}G(\lambda,s)$, $(s_1,s)\in[0,t]^2$,   such that
\begin{equation}\label{DG(s)}
\mathds{E}\left[\int\limits_{[0,t]^{\scriptscriptstyle{2}}}(D_{s_1}G(\lambda,s))^2ds_1ds\right]<+\infty.
\end{equation}
We first state the following Lemma.

\begin{lemma}\label{Lemma(A4)}

 For every $q>1$, there exists a constant $C_{\scriptscriptstyle{0}}(q)$  such that
\begin{equation}\label{etoile(1)}
\sup_{0< r\leq t\leq
1,\lambda\in]0,1]}\mathds{E}\left[(\mathcal{E}(\lambda;r,t))^{q}\right]\leq
C_{\scriptscriptstyle{0}}(q).
\end{equation}
\end{lemma}

\begin{proof}[Proof of Lemma \ref{Lemma(A4)}]

Let $\lambda \in]0,1]$ and $q>1$. For fixed  $0< r \leq t \leq 1$,
\eqref{Prop_UpperBound_Dens.(1bis)}  gives
\begin{equation*}
    \mathcal{E}^q(\lambda;r,t)=M(\lambda;r,t,q)\exp\left(\frac{\displaystyle{\lambda^{\scriptscriptstyle{2}}(q^2-q)}}{\displaystyle{2}}\int\limits_r^t(\sigma')^2(\rho\lambda^{\scriptscriptstyle{2}}, \lambda Z_{\rho}+x_{\scriptscriptstyle{0}})d{\rho}+q\lambda^{\scriptscriptstyle{2}}\int\limits_r^tb'(\rho\lambda^{\scriptscriptstyle{2}},
    \lambda Z_{\rho}+x_{\scriptscriptstyle{0}})d{\rho}\right),
\end{equation*}
where
\begin{equation}\label{M(r,s,q)}
M(\lambda;r,t,q)=\exp\left(\int\limits_r^t\lambda q\sigma'(\rho\lambda^{\scriptscriptstyle{2}},
\lambda Z_{\rho}+x_{\scriptscriptstyle{0}})dW_{\rho}-\frac{\displaystyle{1}}{\displaystyle{2}}\int\limits_r^t(q\lambda\sigma')^2(\rho\lambda^{\scriptscriptstyle{2}},
    \lambda Z_{\rho}+x_{\scriptscriptstyle{0}})d{\rho}\right).
\end{equation}
In fact, since $\sigma'$  is bounded, the stochastic exponential $M(\lambda;r,t,q)$ verifies Novikov's condition; therefore
it is a martingale. So  $\mathds{E}(M(\lambda;r,t,q))=1$. In addition, since $b'$ is also bounded and $\lambda \in]0,1]$, we get $ \mathds{E}\left[(\mathcal{E}(\lambda;r,t))^{q}\right]\leq C_{\scriptscriptstyle{0}}(q)$,
where
$C_{\scriptscriptstyle{0}}(q)=\exp\left(2(q^2-q)\|\sigma'\|^2_{\scriptscriptstyle{\infty}}+2q\|b'\|_{\scriptscriptstyle{\infty}}\right)$. Consequently  \eqref{etoile(1)} is established.
\end{proof}
\begin{proposition}\label{Prop(C)}
There is a constant $\mathcal{C}$ (not depending on $x_{\scriptscriptstyle{0}}$)  such that
\begin{equation*}
    \sup_{(t,\lambda)\in]0,1]^2}\mathds{E}[(<DZ_t,DZ_t>)^{-p}]\leq \mathcal{C}, \ \ \forall p \geq 1.
\end{equation*}
\end{proposition}

\begin{proof}[Proof of Proposition \ref{Prop(C)}]

Let $t\in ]0,1]$ fixed, $\epsilon_{\scriptscriptstyle{0}}= \frac{c_{\scriptscriptstyle{0}}t}{8}$, where $c_{\scriptscriptstyle{0}}$ is a non-degeneracy constant of $\sigma^2$  in the sense of Definition \ref{DefNonDeg}. Consider $\epsilon\in]0,\epsilon_{\scriptscriptstyle{0}}[$, we set $N:=N^{\lambda}=<DZ^{\lambda}_t,DZ^{\lambda}_t>$, where we recall that $<DZ_t,DZ_t>$ appears in \eqref{Prop_UpperBound_Dens.(1)} and
\eqref{Prop_UpperBound_Dens.(1bis)}. According to Proposition
\ref{Prop_E(N{-p})}  we have to evaluate ${P}(N\leq
\epsilon)$. Taking into account the lower bound of $\sigma$  we have
\begin{eqnarray}
  {P}(N\leq
\epsilon) &\leq & {P}\left(\int_0^t dr \mathcal{E}^2(\lambda;r,t)\leq \frac{\epsilon}{c_{\scriptscriptstyle{0}}}\right) \label{Prop_UpperBound_Dens.(2)}\\
&\leq &{P}\left({\left(\int_{t-\frac{4\epsilon}{c_{\scriptscriptstyle{0}}}}^t
dr \mathcal{E}^2(\lambda;r,t)\right)}^{\frac{1}{2}}\leq \sqrt{\frac{\epsilon}{c_{\scriptscriptstyle{0}}}}\right)\nonumber\\
   &\leq & {P}\left({\left(\int_{t-\frac{4\epsilon}{c_{\scriptscriptstyle{0}}}}^t dr\right)}^{\frac{1}{2}}-{\left(\int_{t-\frac{4\epsilon}{c_{\scriptscriptstyle{0}}}}^t  \mathcal{E}^2(\lambda;r,t)dr\right)}^{\frac{1}{2}}\geq \sqrt{\frac{\epsilon}{c_{\scriptscriptstyle{0}}}}\right). \nonumber
\end{eqnarray}
By the inverse triangle inequality of the
$L^2([t-\frac{4\epsilon}{c_{\scriptscriptstyle{0}}},t])$-norm  we
get
\begin{equation*}
{P}(N\leq \epsilon) \leq
{P}\left(\int_{t-\frac{4\epsilon}{c_{\scriptscriptstyle{0}}}}^t
(1-\mathcal{E}(\lambda;r,t))^2 dr\geq
\frac{\epsilon}{c_{\scriptscriptstyle{0}}}\right).
\end{equation*}
Let $p\geq 1$. By Chebyshev's inequality  this is lower than
\begin{equation*}
    {\left(\frac{\displaystyle{c_{\scriptscriptstyle{0}}}}{\displaystyle{\epsilon}}\right)}^{\scriptscriptstyle{p+1}}
    \mathds{E}\left[{\left(\int_{t-\frac{4\epsilon}{c_{\scriptscriptstyle{0}}}}^t
(1-\mathcal{E}(\lambda;r,t))^2 dr\right)}^{\scriptscriptstyle{p+1}}\right].
\end{equation*}
Then, using Jensen's inequality, we get
\begin{equation}\label{Prop_UpperBound_Dens.(3)}
 {P}(N\leq \epsilon) \leq   4^p{\left(\frac{\displaystyle{\epsilon}}{\displaystyle{c_{\scriptscriptstyle{0}}}}\right)}^{-1}
    \int_{t-\frac{4\epsilon}{c_{\scriptscriptstyle{0}}}}^t\mathds{E}\left[{\left(1-\mathcal{E}(\lambda;r,t)\right)}^{\scriptscriptstyle{2(p+1)}}\right]dr.
\end{equation}
Furthermore  \eqref{Prop_UpperBound_Dens.(1bis)}  implies that
$\mathcal{E}(\lambda;r,t)$ solves
\begin{equation*}
    \mathcal{E}(\lambda;r,t)=1+\lambda\int_r^t\mathcal{E}(\lambda;r,s)\sigma'(s\lambda^{\scriptscriptstyle{2}},\lambda Z_s+x_{\scriptscriptstyle{0}})dW_s+\lambda^{\scriptscriptstyle{2}}\int_r^t\mathcal{E}(\lambda;r,s)b'(s\lambda^{\scriptscriptstyle{2}},\lambda Z_s+x_{\scriptscriptstyle{0}})ds.
\end{equation*}
Thus
\begin{align}
  \mathds{E}\left[{\left(\mathcal{E}(\lambda;r,t)-1\right)}^{\scriptscriptstyle{2(p+1)}}\right]
  &\leq2^{\scriptscriptstyle{2(p+1)}}\mathds{E}\left[\left|\lambda\int\limits_r^t\mathcal{E}(\lambda;r,s)\sigma'(s\lambda^{\scriptscriptstyle{2}},\lambda Z_s+x_{\scriptscriptstyle{0}})dW_s\right|^{\scriptscriptstyle{2(p+1)}}\right] \nonumber\\
  &+2^{\scriptscriptstyle{2(p+1)}}\mathds{E}\left[\left|\lambda^{\scriptscriptstyle{2}}\int\limits_r^t\mathcal{E}(\lambda;r,s)
  b'(s\lambda^{\scriptscriptstyle{2}},\lambda Z_s+x_{\scriptscriptstyle{0}})ds\right|^{\scriptscriptstyle{2(p+1)}}\right]. \label{7.19bis}
\end{align}
On one hand, using  Jensen's inequality and $\lambda \in ]0,1]$, we obtain
\begin{equation}\label{PropUpB_Ineq1}
\mathds{E}\left[\left|\lambda^{\scriptscriptstyle{2}}\int\limits_r^t\mathcal{E}(\lambda;r,s)b'(s\lambda^{\scriptscriptstyle{2}},\lambda Z_s+x_{\scriptscriptstyle{0}})ds\right|^{\scriptscriptstyle{2(p+1)}}\right]\leq
\|b'\|_{\scriptscriptstyle{\infty}}
|t-r|^{\scriptscriptstyle{2p+1}}\int\limits_r^t\mathds{E}\left[(\mathcal{E}(\lambda;r,s))^{\scriptscriptstyle{2(p+1)}}\right]ds.
\end{equation}
On the other hand, by Burkholder-Davis-Gundy inequality, we get
\begin{equation*}
   \mathds{E}\left[\left|\lambda\int\limits_r^t\mathcal{E}(\lambda;r,s)\sigma'(s\lambda^{\scriptscriptstyle{2}},\lambda Z_s+x_{\scriptscriptstyle{0}})dW_s\right|^{\scriptscriptstyle{2(p+1)}}\right] \leq  \mathds{E}\left[\left|\lambda^{\scriptscriptstyle{2}}\int\limits_r^t\mathcal{E}^{\scriptscriptstyle{2}}(\lambda;r,s)(\sigma')^2
   (s\lambda^{\scriptscriptstyle{2}},\lambda Z_s+x_{\scriptscriptstyle{0}})ds\right|^{\text{p+1}}\right].
\end{equation*}
Applying again Jensen's inequality gives
\begin{equation}\label{PropUpB_Ineq2}
   \mathds{E}\left[\left|\lambda\int\limits_r^t\mathcal{E}(\lambda;r,s)\sigma'(s\lambda^{\scriptscriptstyle{2}},\lambda Z_s+x_{\scriptscriptstyle{0}})dW_s\right|^{\scriptscriptstyle{2(p+1)}}\right] \leq  \|\sigma'\|_{\scriptscriptstyle{\infty}}
|t-r|^{p}\int\limits_r^t\mathds{E}\left[(\mathcal{E}(\lambda;r,s))^{\scriptscriptstyle{2(p+1)}}\right]ds.
\end{equation}
 Therefore
\eqref{PropUpB_Ineq1}, \eqref{PropUpB_Ineq2} and \eqref{7.19bis} lead to
\begin{equation}\label{PropUpB_Ineq3}
    \mathds{E}\left[{\left(\mathcal{E}(\lambda;r,t)-1\right)}^{\scriptscriptstyle{2(p+1)}}\right]\leq
    C(T,\|\sigma'\|_{\scriptscriptstyle{\infty}},\|b'\|_{\scriptscriptstyle{\infty}})|t-r|^{p}\int\limits_r^t\mathds{E}\left[(\mathcal{E}(\lambda;r,s))^{\scriptscriptstyle{2(p+1)}}\right]ds.
\end{equation}
By Lemma \ref{Lemma(A4)} and \eqref{PropUpB_Ineq3}, there is a constant $C_0(2(p+1))$  such that
\begin{equation}\label{PropUpB_Ineq3Bis}
    \mathds{E}\left[{\left(\mathcal{E}(\lambda;r,t)-1\right)}^{\scriptscriptstyle{2(p+1)}}\right]\leq
    C_1(T,\|\sigma'\|_{\scriptscriptstyle{\infty}},\|b'\|_{\scriptscriptstyle{\infty}})C_0(2(p+1)).
\end{equation}
Then, coming back to \eqref{Prop_UpperBound_Dens.(3)} and using \eqref{PropUpB_Ineq3Bis}, we obtain
\begin{equation}\label{Prop_UpperBound_Dens.(3bis)}
\forall \epsilon \in]0,\epsilon_{\scriptscriptstyle{0}}], \ \ {P}(N\leq \epsilon)\leq
C(p)\epsilon^{\scriptscriptstyle{p+1}},
\end{equation}
where
$C(p)=\frac{\displaystyle{4^{\scriptscriptstyle{2(p+1)}}C_{\scriptscriptstyle{0}}(2(p+1))C_1(T,\|\sigma'\|_{\scriptscriptstyle{\infty}},\|b'\|_{\scriptscriptstyle{\infty}})}}
{\displaystyle{p+1}}$. Finally, using Proposition \ref{Prop_E(N{-p})}, the result follows.
\end{proof}

We go on with the proof of Proposition \ref{Prop_UpperBound_Dens} taking into account the considerations before Lemma \ref{Lemma(A4)}. In fact   \cite[Proposition 2.1.1]{LivreNual} allows us to express,  for
fixed $t \in ]0,1]$,
\begin{equation}\label{Prop_UpperBound_Dens.(3bisbis)}
    q_t^{\lambda}(z)=\mathds{E}\left[\mathds{1}_{\{Z_t>z\}}\delta(G(\lambda,\cdot))\right];
\end{equation}
using Cauchy-Schwarz inequality, it implies that
\begin{equation}\label{Prop_UpperBound_Dens.(4)}
     q_t^{\lambda}(z)\leq \sqrt{\mathds{E}\left[\left|\delta(G(\lambda,\cdot))\right|^2\right]}.
\end{equation}
According to (1.48) in \cite{LivreNual},  \eqref{Prop_UpperBound_Dens.(4)} implies
\begin{equation}\label{Prop_UpperBound_Dens.(4bis)}
     q_t^{\lambda}(z)\leq {\left(\mathds{E}\left[\int\limits_0^tG^2(\lambda,s)ds\right]+
    \mathds{E}\left[\int\limits_{[0,t]^{\scriptscriptstyle{2}}}(D_{s_1}G(\lambda,s))^2ds_1ds\right]\right)}^{\scriptscriptstyle{\frac{1}{2}}}.
\end{equation}
Now  we state a result that estimates the two terms in
the right-hand side of \eqref{Prop_UpperBound_Dens.(4bis)}. Indeed,
we have the following.
\begin{proposition}\label{PropEstimI_1.I_2}
For every $\lambda\in]0,1]$, $G(\lambda,\cdot)\in \mathds{L}^{1,2}$. Moreover, the following statements hold:
(i) $\mathds{E}\left[\int\limits_0^tG^2(\lambda,s)ds\right]\leq \mathcal{C}_1\left(1+ |x_{\scriptscriptstyle
   0}|^2\right),\ \ $   (ii) $\mathds{E}\left[\int\limits_{[0,t]^{\scriptscriptstyle{2}}}(D_{s_1}G(\lambda,s))^2ds_1ds\right]\leq \mathcal{C}_2 \left(1+|x_{\scriptscriptstyle 0}|^8\right)$,
%
where  $\mathcal{C}_1$ and $\mathcal{C}_2$ depend on $T$,
$\|\sigma'\|_{\scriptscriptstyle{\infty}}$ and
$\|b'\|_{\scriptscriptstyle{\infty}}$, but not  on $x_0$.
\end{proposition}
\begin{proof}[Proof of Proposition \ref{PropEstimI_1.I_2}]
(i) First  we set
$I_{\scriptscriptstyle{1}}=\mathds{E}\left[\int\limits_0^tG^2(\lambda,s)ds\right]$,  recalling  that
\begin{equation*}
G(\lambda,s)=\frac{\displaystyle{\sigma(s\lambda^{\scriptscriptstyle{2}},\lambda Z_s+x_{\scriptscriptstyle{0}})\mathcal{E}(\lambda;s,t)}}{\displaystyle{G_{den}}},\
\  \mbox{where }\ G_{den}=<DZ_t,DZ_t>.
\end{equation*}
By Cauchy-Schwarz inequality, we have
\begin{equation}\label{I_1}
I_{\scriptscriptstyle{1}}\leq {\left( \mathds{E}\left[\int\limits_0^t \sigma^{\scriptscriptstyle{4}}(s\lambda^{\scriptscriptstyle{2}},\lambda Z_s+x_{\scriptscriptstyle{0}})ds\right] \mathds{E}\left[\int\limits_0^t\frac{\displaystyle{\mathcal{E}^{\scriptscriptstyle{4}}(\lambda;s,t)}}{\displaystyle{G^4_{den}}}
 ds\right]\right)}^{\frac{1}{2}}.
\end{equation}
Since $\sigma$ has linear growth, by Lemma \ref{Lemma(A3)} and using again Cauchy-Schwarz inequality,  there is a constant $C_{\scriptscriptstyle{1}}$ such that
\begin{equation}\label{PropEstimI_1.I_2.(1bis)}
I_{\scriptscriptstyle{1}}\leq C_{\scriptscriptstyle{1}} \left(1+|x_{\scriptscriptstyle
    0}|^2\right){\left(\mathds{E}\left[{G_{den}}^{\scriptscriptstyle{-8}}\right]
\int\limits_0^t
\mathds{E}\left[\mathcal{E}^{\scriptscriptstyle{8}}(\lambda;s,t)
\right]ds\right)}^{\scriptscriptstyle{\frac{1}{4}}}.
\end{equation}
Consequently, using Proposition \ref{Prop(C)} and Lemma \ref{Lemma(A4)}, the first item of Proposition
\ref{PropEstimI_1.I_2} is established.\\

(ii) We set
$I_{\scriptscriptstyle{2}}=\mathds{E}\left[\int\limits_{[0,t]^{\scriptscriptstyle{2}}}(D_{s_1}G(\lambda,s))^2ds_1ds\right]$.

On one hand, by usual Malliavin differentiation rules, we obtain
\begin{align}
D_{s_1}G(\lambda,s)&=&\lambda\sigma'(s\lambda^{\scriptscriptstyle{2}},\lambda Z_s+x_{\scriptscriptstyle{0}})
\sigma(s_{\scriptscriptstyle{1}}\lambda^{\scriptscriptstyle{2}},\lambda Z_{s_1}+x_{\scriptscriptstyle{0}})\mathcal{E}(\lambda;s_{\scriptscriptstyle{1}},s)\frac{\displaystyle{\mathcal{E}(\lambda;s,t)}}
{\displaystyle{G_{den}}} \mathds{1}_{[0,s]}(s_{\scriptscriptstyle{1}})\nonumber \\
& &\label{PropEstimI_1.I_2.(2)} \\
&+&\sigma(s\lambda^{\scriptscriptstyle{2}},\lambda Z_s+x_{\scriptscriptstyle{0}})\frac{\displaystyle{D_{s_1}\mathcal{E}(\lambda;s,t)}}
{\displaystyle{G_{den}}}-\sigma(s\lambda^{\scriptscriptstyle{2}},\lambda Z_s+x_{\scriptscriptstyle{0}})\mathcal{E}(\lambda;s,t)\frac{\displaystyle{D_{s_1}
G_{den}}}{\displaystyle{G^2_{den}}}.\nonumber
\end{align}
The right-hand side being measurable with respect to $\Omega\times[0,T]^2$, $G(\lambda,\cdot)$ will belong to $\mathds{L}^{1,2}$ if (ii) is established. At this point we need to evaluate $D_{s_1}\mathcal{E}(\lambda;s,t)$. From now on, for
the sake of simplicity,  we will only expose the calculations in the
case when $b\equiv 0$. In fact, we have
\begin{align}
  D_{s_1}\mathcal{E}(\lambda;s,t)& = \mathcal{E}(\lambda;s,t) D_{s_1}\left(\lambda\int\limits_s^t\sigma'({\ell}\lambda^{\scriptscriptstyle{2}},\lambda Z_{\ell}+x_{\scriptscriptstyle{0}})dW_{\ell}-\frac{\lambda^{\scriptscriptstyle{2}}}{2}\int\limits_s^t
  (\sigma')^{\scriptscriptstyle{2}}({\ell}\lambda^{\scriptscriptstyle{2}},\lambda Z_{\ell}+x_{\scriptscriptstyle{0}})d\ell\right)\nonumber\\
   &= \mathcal{E}(\lambda;s,t)\left( \vphantom{\int\limits_0^t}\lambda\sigma'(s_{\scriptscriptstyle{1}}\lambda^{\scriptscriptstyle{2}},\lambda Z_{s_1}+x_{\scriptscriptstyle{0}})\mathds{1}_{[s,t]}(s_{\scriptscriptstyle{1}})\right.
   \nonumber \\
   &+ \lambda^{\scriptscriptstyle{2}}\sigma(s_{\scriptscriptstyle{1}}\lambda^{\scriptscriptstyle{2}},\lambda Z_{s_1}+x_{\scriptscriptstyle{0}})\int\limits_s^t\mathds{1}_{[s_{\scriptscriptstyle{1}},t]}(\ell)
   \sigma''({\ell}\lambda^{\scriptscriptstyle{2}},\lambda Z_{\ell}+x_{\scriptscriptstyle{0}})\mathcal{E}(\lambda;s_{\scriptscriptstyle{1}},\ell)dW_{\ell}  \label{PropEstimI_1.I_2.(2bis)}\\
  &-\left.\lambda^{\scriptscriptstyle{3}}
  \sigma(s_{\scriptscriptstyle{1}}\lambda^{\scriptscriptstyle{2}},\lambda Z_{s_1}+x_{\scriptscriptstyle{0}})\int\limits_s^t\mathds{1}_{[s_{\scriptscriptstyle{1}},t]}(\ell)
  (\sigma'\sigma'')({\ell}\lambda^{\scriptscriptstyle{2}},\lambda Z_{\ell}+x_{\scriptscriptstyle{0}})\mathcal{E}(\lambda;s_{\scriptscriptstyle{1}},\ell)d\ell\right)\nonumber.
\end{align}
On the other hand, we get
\begin{align}
D_{s_1}G_{den}&=2\lambda\int\limits_0^t
\sigma\sigma'(\xi\lambda^{\scriptscriptstyle{2}},\lambda Z_{\xi}+x_{\scriptscriptstyle{0}})\mathds{1}_{[s_{\scriptscriptstyle{1}},t]}(\xi)\sigma(s_{\scriptscriptstyle{1}}\lambda^{\scriptscriptstyle{2}},\lambda Z_{s_{\scriptscriptstyle{1}}}+x_{\scriptscriptstyle{0}})\mathcal{E}(\lambda;s_{\scriptscriptstyle{1}},\xi)
\mathcal{E}^{\scriptscriptstyle{2}}(\lambda;\xi,t)d\xi \nonumber\\
&+2\int\limits_0^t
\sigma^{\scriptscriptstyle{2}}(\xi\lambda^{\scriptscriptstyle{2}},\lambda Z_{\xi}+x_{\scriptscriptstyle{0}})\mathcal{E}(\lambda;\xi,t)D_{s_1}\mathcal{E}(\lambda;\xi,t)d\xi. \label{DGden}
\end{align}
Therefore, coming back to \eqref{PropEstimI_1.I_2.(2)} and  using
\eqref{PropEstimI_1.I_2.(2bis)},  we
obtain that
\begin{equation}\label{PropEstimI_1.I_2.(4)}
I_{\scriptscriptstyle{2}}\leq 4\left[J_1+J_2+J_3\right],
\end{equation}
where
\begin{align*}
  J_1 &= \mathds{E}\left[\int\limits_0^tds\int\limits_0^sds_{\scriptscriptstyle{1}}\left|\sigma'(s\lambda^{\scriptscriptstyle{2}},\lambda Z_s+x_{\scriptscriptstyle{0}})\sigma(s_{\scriptscriptstyle{1}}\lambda^{\scriptscriptstyle{2}},\lambda Z_{s_1}+x_{\scriptscriptstyle{0}})\mathcal{E}(\lambda;s_{\scriptscriptstyle{1}},s)
  \frac{\displaystyle{\mathcal{E}(\lambda;s,t)}}
{\displaystyle{G_{den}}} \right|^{{{2}}}\right],\\
J_2 &= \mathds{E}\left[\int\limits_{[0,t]^{\scriptscriptstyle{2}}}ds_{\scriptscriptstyle{1}}ds \sigma^{\scriptscriptstyle{2}}(s\lambda^{\scriptscriptstyle{2}},\lambda Z_s+x_{\scriptscriptstyle{0}})\frac{\displaystyle{\mathcal{E}^{\scriptscriptstyle{2}}(\lambda;s,t)}}
{\displaystyle{G^{\scriptscriptstyle{2}}_{den}}}\left( \vphantom{\int\limits_0^t}\mathds{1}_{[s,t]}(s_{\scriptscriptstyle{1}})
\sigma'(s_{\scriptscriptstyle{1}}\lambda^{\scriptscriptstyle{2}},\lambda Z_{s_1}+x_{\scriptscriptstyle{0}})\right.\right.\\
 &+\left.\left.
\sigma(s_{\scriptscriptstyle{1}}\lambda^{\scriptscriptstyle{2}},\lambda Z_{s_1}+x_{\scriptscriptstyle{0}})\int\limits_s^t\mathds{1}_{[s_{\scriptscriptstyle{1}},t]}(\ell)\sigma''(\ell\lambda^{\scriptscriptstyle{2}},\lambda Z_{\ell}+x_{\scriptscriptstyle{0}})\mathcal{E}(\lambda;s_{\scriptscriptstyle{1}},\ell)
\left[\vphantom{\frac{\displaystyle{\mathcal{E}^{\scriptscriptstyle{2}}(\lambda;s,t)}}
{\displaystyle{G^{\scriptscriptstyle{2}}_{den}}}}dW_{\ell}-\sigma'\sigma''(\ell\lambda^{\scriptscriptstyle{2}},\lambda Z_{\ell}+x_{\scriptscriptstyle{0}})d\ell\right]\right)^{{2}} \vphantom{\frac{\displaystyle{\mathcal{E}^{\scriptscriptstyle{2}}(\lambda;s,t)}}
{\displaystyle{G^{\scriptscriptstyle{2}}_{den}}}}\right],\\
  J_3 &=\mathds{E}\left[\int\limits_{[0,t]^{\scriptscriptstyle{2}}}ds_{\scriptscriptstyle{1}}ds \frac{\displaystyle{\sigma^{\scriptscriptstyle{2}}(s\lambda^{\scriptscriptstyle{2}},\lambda Z_s+x_{\scriptscriptstyle{0}})}}
{\displaystyle{G^4_{den}}}\mathcal{E}^{\scriptscriptstyle{2}}(\lambda;s,t)\left(D_{s_1}G_{den}\right)^2
\right],
\end{align*}
with  $D_{s_1}G_{den}$ given in \eqref{DGden}. In the sequel  we will enumerate constants $K_{\scriptscriptstyle{1}}$ to $K_{\scriptscriptstyle{20}}$; all those will not depend on $x_{\scriptscriptstyle{0}}$ or $t$, but eventually on $T$, $\sigma$ and $b$. We start  estimating $J_1$. Since $\sigma'$ is bounded, by Cauchy-Schwarz inequality,  we have
\begin{equation*}
    J_1\leq
    K_{\scriptscriptstyle{1}}{\left( \mathds{E}\left[ \int\limits_0^t ds\int\limits_0^s ds_{\scriptscriptstyle{1}}\sigma^{\scriptscriptstyle{4}}(s_{\scriptscriptstyle{1}}\lambda^{\scriptscriptstyle{2}},\lambda Z_{s_1}+x_{\scriptscriptstyle{0}}) \right]\mathds{E}\left[\int\limits_0^t ds \int\limits_0^s ds_{\scriptscriptstyle{1}} \mathcal{E}^{\scriptscriptstyle{4}}(\lambda;s,t) \frac{\displaystyle{\mathcal{E}^{\scriptscriptstyle{4}}(\lambda;s_{\scriptscriptstyle{1}},s)}}{\displaystyle{G^4_{den}}}\right]\right)}^{\frac{1}{2}}.
\end{equation*}
Since $\sigma$ has linear growth, Lemma \ref{Lemma(A3)} and a further use of  Cauchy-Schwarz inequality imply that,  $J_1$ is bounded by
\begin{equation*}
    K_{\scriptscriptstyle{2}}\left(1+|x_{\scriptscriptstyle{0}}|^2\right)
    {\left(\mathds{E}\left[G_{den}^{\scriptscriptstyle{-8}} \right]\right)}^{\frac{1}{4}}
    {\left(\mathds{E}\left[\int\limits_0^t ds \int\limits_0^s ds_{\scriptscriptstyle{1}} \mathcal{E}^{\scriptscriptstyle{16}}(\lambda;s,t) \right]\mathds{E}\left[\int\limits_0^t ds\int\limits_0^s ds_{\scriptscriptstyle{1}}  \mathcal{E}^{\scriptscriptstyle{16}}(\lambda;s_{\scriptscriptstyle{1}},s) \right]\right)}^{\frac{1}{8}}.
\end{equation*}
Therefore, by Proposition \ref{Prop(C)} and Lemma \ref{Lemma(A4)}, we obtain
that
\begin{equation}\label{UpBoundJ1}
    J_1\leq K_{\scriptscriptstyle{3}}\left(1+|x_{\scriptscriptstyle{0}}|^2\right).
\end{equation}
We go on with the analysis of $J_2$. Since $\sigma'$,
$\sigma''$ are bounded, we have
\begin{equation*}
J_2\leq
K_{\scriptscriptstyle{4}}\mathds{E}\left[\int\limits_{[0,t]^2}ds ds_{\scriptscriptstyle{1}}\sigma^{\scriptscriptstyle{2}}(s\lambda^{\scriptscriptstyle{2}},\lambda Z_s+x_{\scriptscriptstyle{0}})\frac{\displaystyle{\mathcal{E}^{\scriptscriptstyle{2}}(\lambda;s,t)}}
{\displaystyle{G^{\scriptscriptstyle{2}}_{den}}}\left(\vphantom{\frac{\displaystyle{\mathcal{E}^{\scriptscriptstyle{2}}(\lambda;s,t)}}
{\displaystyle{G^{\scriptscriptstyle{2}}_{den}}}}1+\sigma^{\scriptscriptstyle{2}}(s_{\scriptscriptstyle{1}}\lambda^{\scriptscriptstyle{2}},\lambda Z_{s_{\scriptscriptstyle{1}}}+x_{\scriptscriptstyle{0}})M^{\scriptscriptstyle{2}}(s,s_{\scriptscriptstyle{1}};t)\vphantom{\frac{\displaystyle{\mathcal{E}^{\scriptscriptstyle{2}}(\lambda;s,t)}}
{\displaystyle{G^{\scriptscriptstyle{2}}_{den}}}}\right)\right],
\end{equation*}
where  $M(s,
s_{\scriptscriptstyle{1}};t)=\int\limits_{s\vee
s_{\scriptscriptstyle{1}}}^t\sigma''(\ell\lambda^{\scriptscriptstyle{2}},\lambda Z_{\ell}+x_{\scriptscriptstyle{0}})\mathcal{E}(\lambda;s_{\scriptscriptstyle{1}},\ell)dW_{\ell}$, $\scriptstyle{t\geq s\vee
s_{\scriptscriptstyle{1}}}$, is a  martingale having  all moments because of Lemma \ref{Lemma(A4)}. Since $\sigma$ has linear growth, using Cauchy-Schwarz inequality and Lemma \ref{Lemma(A3)}, we get
\begin{equation*}
    J_2\leq
K_{\scriptscriptstyle{5}}(1+|x_{\scriptscriptstyle{0}}|^4)\left(\mathds{E}\left[\int\limits_0^tds\mathcal{E}^{\scriptscriptstyle{16}}(\lambda;s,t)\right]
\mathds{E}\left[G_{den}^{\scriptscriptstyle{-16}}\right]\right)^{\frac{1}{8}}\left(\int\limits_{[0,t]^2}ds ds_{\scriptscriptstyle{1}}E[M^{\scriptscriptstyle{8}}(s,
s_{\scriptscriptstyle{1}};t)]\right)^{\frac{1}{4}}.
\end{equation*}
Then   Burkholder-Davis-Gundy inequality, Lemma \ref{Lemma(A4)} and Proposition \ref{Prop(C)} imply
\begin{equation}\label{J2UpBound}
    J_2\leq K_{\scriptscriptstyle{6}}(1+|x_{\scriptscriptstyle{0}}|^4).
\end{equation}
Finally  we treat $J_3$. Applying  Cauchy-Schwarz  inequality, we have
\begin{equation*}
    J_3\leq \left(\mathds{E}\left[\int\limits_{[0,t]^2}ds ds_{\scriptscriptstyle{1}}\frac{\sigma^{\scriptscriptstyle{4}}(s\lambda^{\scriptscriptstyle{2}},\lambda Z_s+x_{\scriptscriptstyle{0}})}{G_{den}^{\scriptscriptstyle{8}}}
    \mathcal{E}^{\scriptscriptstyle{4}}(\lambda;s,t)\right]\mathds{E}\left[\int\limits_{[0,t]^2}ds ds_{\scriptscriptstyle{1}}
(D_{s_1}G_{den})^{\scriptscriptstyle{4}}\right]\right)^{\frac{1}{2}}.
\end{equation*}
Since $\sigma$ has linear growth, again by Cauchy-Schwarz inequality and Lemma \ref{Lemma(A3)}, we get
\begin{equation*}
    J_3\leq K_{\scriptscriptstyle{7}}(1+|x_{\scriptscriptstyle{0}}|^2)\left(\mathds{E}\left[G_{den}^{\scriptscriptstyle{-32}}\right]
    \int\limits_0^tds \mathds{E}\left[\mathcal{E}^{\scriptscriptstyle{16}}(\lambda;s,t)\right]\right)^{\frac{1}{8}}
    \left(\mathds{E}\left[\int\limits_0^t ds_{\scriptscriptstyle{1}}
(D_{s_1}G_{den})^{\scriptscriptstyle{4}}\right]\right)^{\frac{1}{2}};
\end{equation*}
 by Lemma \ref{Lemma(A4)} and Proposition \ref{Prop(C)}  it follows
\begin{equation}\label{UpBoundIntermJ3}
    J_3\leq K_{\scriptscriptstyle{8}}(1+|x_{\scriptscriptstyle{0}}|^2)\left(\mathds{E}\left[\int\limits_0^t ds_{\scriptscriptstyle{1}}
(D_{s_1}G_{den})^{\scriptscriptstyle{4}}\right]\right)^{\frac{1}{2}}.
\end{equation}
Since $\sigma'$ is bounded, \eqref{DGden} and Jensen's inequality give
\begin{equation}\label{UpBoundDGden4}
    \mathds{E}\left[\int\limits_0^t ds_{\scriptscriptstyle{1}}
(D_{s_1}G_{den})^{\scriptscriptstyle{4}}\right]\leq K_{\scriptscriptstyle{9}}\left(A_1+A_2\right),
\end{equation}
where
\begin{align*}
A_1&= \mathds{E}\left[\int\limits_0^t ds_{\scriptscriptstyle{1}}\int\limits_{s_1}^t d\xi \sigma^{\scriptscriptstyle{4}}(s_{\scriptscriptstyle{1}}\lambda^{\scriptscriptstyle{2}},\lambda Z_{s_{\scriptscriptstyle{1}}}+x_{\scriptscriptstyle{0}}) \sigma^{\scriptscriptstyle{4}}(\xi\lambda^{\scriptscriptstyle{2}},\lambda Z_{\xi}+x_{\scriptscriptstyle{0}})\mathcal{E}^{\scriptscriptstyle{4}}(\lambda;s_{\scriptscriptstyle{1}},\xi)
\mathcal{E}^{\scriptscriptstyle{8}}(\lambda;\xi,t)\right],\\
A_2&=\mathds{E}\left[\int\limits_0^t ds_{\scriptscriptstyle{1}}\int\limits_{0}^t d\xi \sigma^{\scriptscriptstyle{8}}(\xi\lambda^{\scriptscriptstyle{2}},\lambda Z_{\xi}+x_{\scriptscriptstyle{0}})
\mathcal{E}^{\scriptscriptstyle{4}}(\lambda;\xi,t)\left(D_{s_1}\mathcal{E}(\lambda;\xi,t)\right)^{\scriptscriptstyle{4}}\right].
\end{align*}
Since $\sigma$ has linear growth, Cauchy-Schwarz inequality and Lemma \ref{Lemma(A3)} imply that $A_1$ is bounded by
\begin{equation*}
     K_{\scriptscriptstyle{10}}\left(1+|x_{\scriptscriptstyle{0}}|^8\right)\left(\mathds{E}\left[\int\limits_0^t ds_{\scriptscriptstyle{1}}\int\limits_{s_1}^t d\xi \mathcal{E}^{\scriptscriptstyle{8}}(\lambda;s_{\scriptscriptstyle{1}},\xi) \mathcal{E}^{\scriptscriptstyle{16}}(\lambda;\xi,t)\right]\right)^{\frac{1}{2}}.
\end{equation*}
Again, by Cauchy-Schwarz inequality and Lemma \ref{Lemma(A4)}, we obtain
\begin{equation}\label{UpBoundA1}
    A_1\leq K_{\scriptscriptstyle{11}}\left(1+|x_{\scriptscriptstyle{0}}|^8\right).
\end{equation}
We proceed  estimating $A_2$. Using Cauchy-Schwarz inequality, $A_2$ is bounded by
\begin{equation*}
    K_{\scriptscriptstyle{12}}\left(\mathds{E}\left[\int\limits_{[0,t]^2}ds_{\scriptscriptstyle{1}}d\xi \sigma^{\scriptscriptstyle{16}}(\xi\lambda^{\scriptscriptstyle{2}},\lambda Z_{\xi}+x_{\scriptscriptstyle{0}})\right]
    \mathds{E}\left[\int\limits_{[0,t]^2}ds_{\scriptscriptstyle{1}}d\xi \mathcal{E}^{\scriptscriptstyle{8}}(\lambda;\xi,t)\left(D_{s_1}\mathcal{E}(\lambda;\xi,t)\right)^{\scriptscriptstyle{8}}\right]\right)^{\frac{1}{2}}.
\end{equation*}
Since $\sigma$ has linear growth, Cauchy-Schwarz inequality and Lemma \ref{Lemma(A3)} lead to
\begin{equation*}
    A_2 \leq K_{\scriptscriptstyle{13}}\left(1+|x_{\scriptscriptstyle{0}}|^8\right)\left( \mathds{E}\left[\int\limits_{[0,t]^2}ds_{\scriptscriptstyle{1}}d\xi \mathcal{E}^{\scriptscriptstyle{16}}(\lambda;\xi,t)\right] \mathds{E}\left[\int\limits_{[0,t]^2}ds_{\scriptscriptstyle{1}}d\xi \left(D_{s_1}\mathcal{E}(\lambda;\xi,t)\right)^{\scriptscriptstyle{16}}\right]  \right)^{\frac{1}{4}};
\end{equation*}
Lemma \ref{Lemma(A4)} implies
\begin{equation}\label{UpBoundIntermA2}
     A_2 \leq K_{\scriptscriptstyle{14}}\left(1+|x_{\scriptscriptstyle{0}}|^8\right)\left(\mathds{E}\left[\int\limits_{[0,t]^2}ds_{\scriptscriptstyle{1}}d\xi \left(D_{s_1}\mathcal{E}(\lambda;\xi,t)\right)^{\scriptscriptstyle{16}}\right]  \right)^{\frac{1}{4}}.
\end{equation}
Since $\sigma'$ and $\sigma''$ are bounded, using  \eqref{PropEstimI_1.I_2.(2bis)} and Jensen's inequality, it follows that
\begin{equation}\label{DDol}
    \mathds{E}\left[\int\limits_{[0,t]^2}ds_{\scriptscriptstyle{1}}d\xi \left(D_{s_1}\mathcal{E}(\lambda;\xi,t)\right)^{\scriptscriptstyle{16}}\right]\leq  K_{\scriptscriptstyle{15}}\left(R_1+R_2\right),
\end{equation}
where
\begin{align*}
    R_1&= \mathds{E}\left[\int\limits_{[0,t]^2}ds_{\scriptscriptstyle{1}}d\xi \mathcal{E}^{\scriptscriptstyle{16}}(\lambda;\xi,t)\right],\\
    R_2&=\mathds{E}\left[\int\limits_{[0,t]^2}ds_{\scriptscriptstyle{1}}d\xi  \sigma^{\scriptscriptstyle{16}}(s_{\scriptscriptstyle{1}}\lambda^{\scriptscriptstyle{2}},\lambda Z_{s_{\scriptscriptstyle{1}}}+x_{\scriptscriptstyle{0}})\mathcal{E}^{\scriptscriptstyle{16}}(\lambda;\xi,t)
    \left(M^{\scriptscriptstyle{16}}(s_{\scriptscriptstyle{1}},\xi;t)+\int\limits_{\xi\vee s_1}^t\mathcal{E}^{\scriptscriptstyle{16}}(\lambda;s_{\scriptscriptstyle{1}},\rho)d\rho\right)\right],
\end{align*}
and $M(s_{\scriptscriptstyle{1}},\xi;t)=\int\limits_{\xi\vee s_1}^t\sigma''(\rho\lambda^{\scriptscriptstyle{2}},\lambda Z_{\rho}+x_{\scriptscriptstyle{0}})\mathcal{E}(\lambda;s_{\scriptscriptstyle{1}},\rho)dW_{\rho}$, $\scriptstyle{t\geq \xi\vee s_1}$, is a  square integrable martingale,  taking into account   Lemma \ref{Lemma(A4)}. Again by Lemma \ref{Lemma(A4)}, $R_1$ is uniformly bounded in $t$ and $x_{\scriptscriptstyle{0}}$. On the other hand, using Cauchy-Schwarz and Jensen's inequalities,  $R_2$ is bounded by
\begin{equation*}
   K_{\scriptscriptstyle{16}} \left(\mathds{E}\left[\int\limits_{[0,t]^2}ds_{\scriptscriptstyle{1}}d\xi  \sigma^{\scriptscriptstyle{32}}(s_{\scriptscriptstyle{1}}\lambda^{\scriptscriptstyle{2}},\lambda Z_{s_{\scriptscriptstyle{1}}}+x_{\scriptscriptstyle{0}})\mathcal{E}^{\scriptscriptstyle{32}}(\lambda;\xi,t)\right]
   \mathds{E}\left[\int\limits_{[0,t]^2}ds_{\scriptscriptstyle{1}}d\xi \left(M^{\scriptscriptstyle{32}}(s_{\scriptscriptstyle{1}},\xi;t)+\int\limits^t_{\xi\vee s_1}\mathcal{E}^{\scriptscriptstyle{32}}(\lambda;s_{\scriptscriptstyle{1}},\rho)d\rho\right)\right]\right)^{\frac{1}{2}}.
\end{equation*}
Since $\sigma$ has linear growth, again by  Cauchy-Schwarz inequality,  Lemma \ref{Lemma(A3)} and Lemma \ref{Lemma(A4)},  we get
\begin{equation*}
    R_2\leq  K_{\scriptscriptstyle{17}} \left(1+|x_{\scriptscriptstyle{0}}|^{16}\right)\mathds{E}\left[\int\limits_{[0,t]^2}ds_{\scriptscriptstyle{1}}d\xi M^{\scriptscriptstyle{32}}(s_{\scriptscriptstyle{1}},\xi;t)\right].
\end{equation*}
Burkholder-Davis-Gundy  inequality gives
\begin{equation}\label{UpBoundR2}
    R_2\leq K_{\scriptscriptstyle{18}} \left(1+|x_{\scriptscriptstyle{0}}|^{16}\right).
\end{equation}
Coming back to \eqref{UpBoundIntermA2}, using \eqref{UpBoundR2} and \eqref{DDol}, we obtain
\begin{equation}\label{UpBoundA2}
    A_2\leq  K_{\scriptscriptstyle{19}} \left(1+|x_{\scriptscriptstyle{0}}|^{12}\right);
\end{equation}
thus, replacing \eqref{UpBoundA1} and \eqref{UpBoundA2} in \eqref{UpBoundDGden4} and coming back to \eqref{UpBoundIntermJ3}, imply
\begin{equation}\label{J3UpBound}
    J_3\leq K_{\scriptscriptstyle{20}} \left(1+|x_{\scriptscriptstyle{0}}|^{8}\right).
\end{equation}
Consequently, substituting \eqref{UpBoundJ1}, \eqref{J2UpBound} and \eqref{J3UpBound} in \eqref{PropEstimI_1.I_2.(4)}, item (ii) of Proposition \ref{PropEstimI_1.I_2} is established.
\end{proof}
Returning to the proof of Proposition
\ref{Prop_UpperBound_Dens} and substituting in
\eqref{Prop_UpperBound_Dens.(4bis)} the right-hand side of  the first and the second item
of Proposition \ref{PropEstimI_1.I_2},  the inequality
\eqref{DensUpperBound} is  verified. Finally  this  concludes the proof
of Proposition \ref{Prop_UpperBound_Dens}.
\end{proof}



\bibliographystyle{amsplain}
\bibliography{BibFDE}


%
%
%
%
%
%
%


\ACKNO{The  authors were partially supported
by the ANR Project MASTERIE 2010 BLAN 0121 01.   Part of the work was  done during the stay  of the second named author  at Bielefeld University, SFB $701$ (Mathematik). The  authors
acknowledge the stimulating remarks of  an anonymous Referee and
of the Editors. They are also grateful to Dr. Juliet Ryan for her precious
help in correcting several language mistakes.}


\end{document}